\documentclass[preprint,12pt, nonatbib]{elsarticle}
\usepackage{mathscinet}
\usepackage{amssymb}
\usepackage{amsmath}
\usepackage{amsthm}
\usepackage{esint}
\usepackage{mathtools}
\usepackage{bm}
\usepackage[colorlinks]{hyperref}
\AtBeginDocument{
  \hypersetup{
    linkcolor=blue,
    citecolor=green,
  }
}

\makeatletter
\let\c@author\relax
\makeatother
\usepackage[backend=bibtex,style=numeric,sorting=nyt,
giveninits=true,isbn=false,doi=false,url=true,maxbibnames=99]{biblatex}
\renewbibmacro{in:}{\ifentrytype{article}{}{\printtext{\bibstring{in}\intitlepunct}}}
\DeclareFieldFormat*{title}{#1}
\usepackage{cleveref}
\usepackage{geometry}
\usepackage[T1]{fontenc}
\usepackage[mathlines,pagewise]{lineno}
\newcommand*\patchAmsMathEnvironmentForLineno[1]{%
  \expandafter\let\csname old#1\expandafter\endcsname\csname #1\endcsname
  \expandafter\let\csname oldend#1\expandafter\endcsname\csname end#1\endcsname
  \renewenvironment{#1}%
     {\linenomath\csname old#1\endcsname}%
     {\csname oldend#1\endcsname\endlinenomath}}%
\newcommand*\patchBothAmsMathEnvironmentsForLineno[1]{%
  \patchAmsMathEnvironmentForLineno{#1}%
  \patchAmsMathEnvironmentForLineno{#1*}}%
\AtBeginDocument{%
\patchBothAmsMathEnvironmentsForLineno{equation}%
\patchBothAmsMathEnvironmentsForLineno{align}%
\patchBothAmsMathEnvironmentsForLineno{flalign}%
\patchBothAmsMathEnvironmentsForLineno{alignat}%
\patchBothAmsMathEnvironmentsForLineno{gather}%
\patchBothAmsMathEnvironmentsForLineno{multline}%
}

\crefname{equation}{}{}

\newtheorem{theorem}{Theorem}[section]
\newtheorem{lemma}[theorem]{Lemma}

\newtheorem{corollary}[theorem]{Corollary}

\theoremstyle{definition}
\newtheorem{definition}[theorem]{Definition}

\newtheorem{notation}[theorem]{Notation}

\theoremstyle{remark}
\newtheorem{remark}[theorem]{Remark}

\numberwithin{equation}{section}

\journal{~~}

\begin{document}

\begin{frontmatter}



\title{Interior pointwise regularity for elliptic and parabolic equations in divergence form and applications to nodal sets \tnoteref{t1}}

\author[rvt]{Yuanyuan Lian\corref{cor1}}
\ead{lianyuanyuan.hthk@gmail.com; yuanyuanlian@correo.ugr.es}
\cortext[cor1]{MR Author ID:
\href{https://mathscinet.ams.org/mathscinet/2006/mathscinet/search/author.html?mrauthid=1378049}
{1378049};
ORCID: \href{https://orcid.org/0000-0002-2276-3063}{0000-0002-2276-3063}
}
\tnotetext[t1]{This research has been financially supported by the Project PID2020-117868GB-I00 funded by MCIN/AEI /10.13039/501100011033.}

\address[rvt]{Departamento de An\'{a}lisis Matem\'{a}tico, Instituto de Matem\'{a}ticas IMAG, Universidad de Granada}

\begin{abstract}
In this paper, we obtain the interior pointwise $C^{k,\alpha}$ ($k\geq 0$, $0<\alpha<1$) regularity for weak solutions of elliptic and parabolic equations in divergence form. The compactness method and perturbation technique are employed. The pointwise regularity is proved in a very simple way and the results are optimal. In addition, these pointwise regularity can be used to characterize the structure of the nodal sets of solutions.
\end{abstract}


\end{frontmatter}

%%\linenumbers

\section{Introduction}
\label{intro}
In this paper, we investigate the interior pointwise regularity for weak solutions of elliptic equations:
\begin{equation}\label{e.div-0}
(a^{ij}u_i+d^ju)_j+b^iu_i+cu=f-f^i_i ~~\mbox{in}~~B_1
\end{equation}
and parabolic equations:
\begin{equation}\label{e.div}
u_t-(a^{ij}u_i+d^ju)_j+b^iu_i+cu=f-f^i_i ~~\mbox{in}~~Q_1,
\end{equation}
where $B_1\subset \mathbb{R}^n$ is the unit ball and $Q_1\subset \mathbb{R}^{n+1}$ is the unit parabolic round cube ($n\geq 2$). In above equations, the Einstein summation convention is used and hereinafter. Since the proofs are similar, we only give the detailed proof for \cref{e.div}. The analogous results for \cref{e.div-0} (see \Cref{sec:4}) can be obtained more easily and we only point out the major difference in the proof.

Throughout this paper, we always assume that the matrix $\pmb{a}$ is uniformly elliptic with ellipticity constants $0<\lambda\leq \Lambda$ and
\begin{equation}\label{e1.2}
 \pmb{b},\pmb{d}\in L^{2p,2p'}(Q_1),~~ c\in L^{p,p'}(Q_1) ,~~f\in L^{q,q'}(Q_1),~~
  \pmb{f}\in L^2(Q_1),
\end{equation}
where $\pmb{b}=(b^1,...,b^n)$, $\pmb{d}=(d^1,...,d^n)$, $\pmb{f}=(f^1,...,f^n)$,
\begin{equation}\label{e1.3}
\frac{n}{p}+\frac{2}{p'}= 2~ \left(p>\frac{n}{2}\right)
\end{equation}
and
\begin{equation}\label{e1.4}
\frac{n}{q}+\frac{2}{q'}=\frac{n}{2}+2~\left(q\geq \frac{2n}{n+2}, q'\geq 1 ~~\mbox{for}~~n\geq 3
~~\mbox{and}~~q>1, q'\geq 1~~\mbox{for}~~n=2\right).
\end{equation}

\begin{remark}\label{re1.0}
The $f\in L^{s,s'}(U)$ ($U\subset \mathbb{R}^{n+1}$) means
\begin{equation*}
\|f\|_{L^{s,s'}(U)}
:=\left(\int_{\mathbb{R}}\|f\cdot \textbf{1}_{U}\|^{s'}_{L^s(\mathbb{R}^n)}dt\right)^{\frac{1}{s'}}
=\left(\int_{\mathbb{R} }\left(\int_{\mathbb{R}^n}
|f(x,t)|^s\cdot \textbf{1}_{U}dx\right)^{\frac{s'}{s}}dt\right)^{\frac{1}{s'}}<\infty,
\end{equation*}
where $ \textbf{1}_U$ denotes the indicator function in $U$. Above notation is used widely for parabolic equations (see \cite{MR0241822,MR0150444}). It is due to the
occurrence of a function having different integrability with respect to $x$ and $t$ when treating parabolic equations (e.g., the weak solution $u$). If $s=s'$, $f\in L^{s,s'}$ is shortened to $f\in L^s$ in this paper.
\end{remark}

\begin{remark}\label{re1.2}
The conditions \crefrange{e1.2}{e1.4} for $\pmb{b},c,\pmb{d}, f$ and $\pmb{f}$ are probable the minimal requirements to establish the energy inequality (see \cite[Chapter III]{MR0241822}), which is the most fundamental estimate for \cref{e.div}. Our regularity results are derived based on the very conditions.
\end{remark}

\begin{remark}\label{re1.1}
The indices $p,p'$ and $q,q'$ are fixed throughout this paper. We point out that the coefficients $\pmb{b},c$ and $\pmb{d}$ can lie in different function spaces such as $\pmb{b}\in L^{p_b,p_b'}(Q_1)$, $c\in L^{p_c,p_c'}(Q_1)$ and $\pmb{d}\in L^{p_d,p_d'}(Q_1)$ respectively as long as the conditions
\begin{equation*}
  \frac{n}{p_b}+\frac{2}{p_b'}= 1,~
  ~\frac{n}{p_c}+\frac{2}{p_c'}= 2,~
  ~\frac{n}{p_d}+\frac{2}{p_d'}= 1
\end{equation*}
are met. A particular case of \cref{e1.3} and \cref{e1.4} is
\begin{equation}\label{e1.5}
 \pmb{b},\pmb{d}\in L^{n+2}(Q_1),~~ c\in L^{\frac{n+2}{2}}(Q_1),~~
f\in L^{\frac{2(n+2)}{n+4}}(Q_1),~~ \pmb{f}\in L^2(Q_1).
\end{equation}

In addition, it is obvious that \cref{e1.3} and \cref{e1.4} can be replaced by
\begin{equation*}
\frac{n}{p}+\frac{2}{p'}\leq 2,~~\frac{n}{q}+\frac{2}{q'}\leq \frac{n}{2}+2.
\end{equation*}
\end{remark}
~\\

The Schauder estimates (i.e., $C^{k,\alpha}$ estimates for $k\geq 0$ and $0<\alpha<1$) play an important role in the theory for linear and nonlinear elliptic and parabolic equations. For equations with constant coefficients, the estimates can be obtained via
\begin{itemize}
  \item Potential theoretic methods based on the fundamental solution (see \cite[Chapter 4]{MR1814364}, \cite[Chapter 4.3]{MR181836}, \cite[Theorem 1.1.4]{MR1613650}, \cite[Section 1]{MR1658901}).
  \item Energy inequality (see \cite[Chapter 3.4]{MR1669352}, \cite{MR1399194}, \cite[Chapter 6.2, Chapter 7.2]{MR2309679}).
  \item Blow-up argument (see \cite{MR1459795}).
  \item Convolution operator (see \cite{MR209917}).
  \item Mollification technique (see \cite{MR912926}).
  \item Perturbation technique (see \cite{MR2273802}).
\end{itemize}
For general equations, the Schauder estimates can be proved by the technique of freezing coefficients (see \cite[Chapter 3.4]{MR1669352}, \cite[Chapter 6]{MR1814364}, \cite[Section 2]{MR1658901}, \cite[Theorem 4.8]{MR1465184}, \cite{MR2309679}, \cite[Chapter 4.4]{MR181836},\cite[Chapter 6.3]{MR2309679}).

There is another powerful technique. The $C^{k,\alpha}$ space can be characterized by polynomials, which is due to Campanato \cite{MR156188, MR167862}. Caffarelli \cite{MR1005611} (see also \cite[Chapter 8]{MR1351007}) and Safonov \cite[Section 7]{MR1655532} first used this polynomial characterization to prove pointwise Schauder estimates. This technique has several advantages:
\begin{itemize}
  \item We can obtain pointwise regularity.
  \item The method is appropriate for weak solutions.
  \item We obtain the any order regularity directly from the solution itself, i.e., the derivatives of solutions do not appear in the proof.
  \item The assumptions on the coefficients, the right-hand term etc. are usually the weakest since the method mainly uses the scaling property of equations.
\end{itemize}
This technique has became a standard method and was applied widely (e.g. see \cite{MR1713596,lian2020pointwise,MR3980853,MR2334822,MR1606359,MR3246039} for elliptic equations; see \cite{LZ_2022,MR3713553,MR1139064,MR1151267} for parabolic equations).

The pointwise regularity shows clearly that the behavior of a solution near a point is essentially determined by the corresponding coefficients and prescribed date near the same point. In addition, the classical local and global regularity can be obtained directly from the pointwise regularity. Moreover, the pointwise regularity has applications to other important areas such as the free boundary problems \cite{MR2511866, MR3158506, MR3385162, MR4644419}\cite[Proof of Theorem 1.2]{MR3198649} and the structure of nodal sets \cite{MR1794573}.

Han \cite{MR1658901,MR1794573} proved a higher order pointwise regularity based on potential analysis for elliptic and parabolic equations in non-divergence form. He showed that if the solution $u$ vanishes at $0$ up to some order, the smoothness requirements on the coefficients can be reduced (see \cite[Theorem 1.2]{MR1794573} and \cite[Theorem 3.1]{MR1658901}). Then he used this pointwise regularity to give a characterization of the structure of the nodal sets for elliptic equations (see \cite[Theorem 5.1]{MR1794573}).

In this paper, we use the compactness method and perturbation technique along with very simple proofs to obtain the pointwise regularity for weak solutions. Similar to \cite{MR1658901,MR1794573}, we obtained that the smooth conditions for the coefficients can be reduced by order $k$ for $\pmb{a},\pmb{b}$ and $k+1$ for $c, \pmb{d}$ if the solution $u$ vanishes at $0$ up to order $k\geq 1$ (see \Cref{t-Cka2} and \Cref{t-Cka2-ell}). Then we use this pointwise regularity to give an characterization of the structure of the nodal sets (see \Cref{th1.3} and \Cref{th1.4}).

 With the perturbation technique, we regard \cref{e.div-0} (resp. \cref{e.div}) as a perturbation of the Laplace equation (resp. the heat equation). We do not use the Green function nor the potential analysis, which are relatively complicated. Since our proofs rely solely on the energy inequality and the scaling argument, the conditions for the coefficients and the right-hand term are probable the weakest (note that $\pmb{b},c\in L^{\infty}$ are assumed in \cite{MR1658901,MR1794573}). This paper has been inspired by \cite{MR1005611}, \cite{MR1139064} and \cite{MR4644419}.

In addition, we give a better characterization (compare \Cref{th1.4} with \cite[Theorem 5.1]{MR1794573}) of structure of the nodal sets by the Whitney extension and the implicit function theorem, which is relatively simple and inspired by its application in free boundary problems (see \cite[Theorem 7.9]{MR2962060}, \cite{MR3385162}). It is worth pointing out that the pointwise regularity is sufficient and necessary to make an Whitney extension (see \Cref{coA.1}). It indicates that the polynomial characterization of $C^{k,\alpha}$ by Campanato is more essential than the classical definition.

Now, we introduce some notions and definitions. We first introduce the definition of the weak solution
(see \cite[P. 136]{MR0241822} and \cite[P. 102]{MR1465184}).
\begin{definition}\label{de1.weak}
Let $U=\Omega\times (T_1,T_2]$. A function $u\in L^2(T_1,T_2;H_0^1(\Omega))\cap L^{\infty}(T_1,T_2;L^2(\Omega))$ is called a weak solution of \cref{e.div} in $U$ if
\begin{equation}\label{e2.1}
-\int_{U}u\varphi_t+\int_{U}a^{ij}u_j\varphi_i+b^iu_i\varphi+cu\varphi+d^iu\varphi_i =\int_{U}f\varphi+f^i\varphi_i
\end{equation}
for any $\varphi\in L^2(T_1,T_2;H_0^1(\Omega)),~\varphi_t\in L^2(U)$ and $\varphi(\cdot,t)\equiv 0$ for $t=T_1$ and $T_2$.
\end{definition}

\begin{remark}\label{re.weak.0}
We use the standard symbol $L^{s'}(T_1,T_2;X)$ for studying parabolic equations (see \cite[Chapter 5.9.2]{MR2597943}). If $X=L^s$, it can be abbreviated as $L^{s,s'}$ (as shown in \Cref{re1.0}). These two symbols are used both in this paper.
\end{remark}

%\begin{remark}\label{re.weak.1}
%By interpolation, $u\in L^{p_0,p'_0}(U)$, where
%\begin{equation*}
%  \frac{n}{p_0}+\frac{2}{p'_0}=\frac{n}{2},\quad 2\leq p_0\leq \frac{2n}{n-2}.
%\end{equation*}
%\end{remark}

\begin{remark}\label{re.weak.2}
There are other definitions of weak solutions (see \cite[the Definition on P. 374]{MR2597943}, \cite[Definition 3.1.1]{MR2309679}). Since the condition on the coefficients is the weakest one among various definitions, we use \Cref{de1.weak} in this paper.
\end{remark}
~\\

Next, we introduce some definitions of pointwise smoothness for a function. First, for $U=\Omega\times (T_1,T_2]$, let
\begin{equation*}
\|f\|^*_{L^{s,s'}(U)}:=\left(\fint^{T_2}_{T_1}  \left(\fint_{\Omega} |f|^s\right)^{s'/s} \right)^{1/s'},
\end{equation*}
where $\fint$ denotes the integral average. The benefit of this notation is that the scaling of $\|f\|^*_{L^{s,s'}(U)}$ is exactly that of $f$, which can make the proof concise when we deal with the scaling argument (see the proofs of \Cref{t-Ca-scaling}, \Cref{t-C1a-scaling} and \Cref{t-Cka-scaling}).

A function is $C^{k,\alpha}$ at one point is as follows.
\begin{definition}\label{d-f1}
Let $U\subset \mathbb{R}^{n+1}$ be a bounded domain, $f:U \rightarrow \mathbb{R}$ be a function and $0<\alpha\leq 1$. We say that $f$ is $C_{s,s'}^{k,\alpha}$ ($k\geq 0$) at $X_0\in U$ or $f\in C_{s,s'}^{k, \alpha}(X_0)$ if there exist $K,r_0>0$ and $P\in \mathcal{P}_k$ ($\mathcal{P}_k$ denotes the set of parabolic polynomials of degree less than or equal to $k$, see \Cref{no1.1}) such that $Q_r(X_0)\subset U$ and
\begin{equation}\label{m-holder}
  \|f-P\|^*_{L^{s,s'}(Q_{r}(X_0))}\leq K r^{k+\alpha},~~\forall~0<r<r_0.
\end{equation}
Then we define for any $1\leq m\leq k$,
\begin{equation*}
  D^mf(X_0)=D^mP(X_0),~~\|f\|_{C_{s,s'}^{k}(X_0)}=\sum_{m=0}^{k}|D^m P(X_0)|,
\end{equation*}

\begin{equation*}
[f]_{C_{s,s'}^{k,\alpha}(X_0)}=\min \left\{K\big | \cref{m-holder} ~\mbox{holds with}~K\right\}
\end{equation*}
and
\begin{equation*}
\|f\|_{C_{s,s'}^{k, \alpha}(X_0)}=\|f\|_{C_{s,s'}^{k}(X_0)}+[f]_{C_{s,s'}^{k, \alpha}(X_0)}.
\end{equation*}
If $f\in C_{s,s'}^{k, \alpha}(X)$ for any $X\in U'\subset U$ with the same $r_0$ and
\begin{equation*}
  \|f\|_{C_{s,s'}^{k,\alpha}(\bar{U}')}:= \sup_{X\in U'} \|f\|_{C_{s,s'}^{k}(X)}+\sup_{X\in U'} [f]_{C_{s,s'}^{k,\alpha}(X)}<+\infty,
\end{equation*}
we say that $f\in C_{s,s'}^{k,\alpha}(\bar{U}')$.
\end{definition}

%\begin{remark}\label{r-df1.1}
%If $U$ is a bounded smooth domain, the definitions of $C_{s,s'}^{k,\alpha}(\bar{U})$ is equivalent to the above usual classical definitions.
%\end{remark}

\begin{remark}\label{r-df1.2}
The above definition $f\in C_{s,s'}^{k,\alpha}(\bar{U}')$ is equivalent to the classical definition. Since we can not find any existing literature to cite, we give the detailed proof in the Appendix (see \Cref{S.A}).
\end{remark}

\begin{remark}\label{r-df1.6}
If $\pmb{f}$ is vector valued function, $\pmb{f}\in C_{s,s'}^{k, \alpha}(X_0)$ means that each component of $\pmb{f}$ is $C_{s,s'}^{k, \alpha}$ at $(X_0)$ and $\|\pmb{f}\|_{C_{s,s'}^{k, \alpha}(X_0)}$ can be defined correspondingly.
\end{remark}

~\\

%\begin{remark}\label{r-df1.2}
%If we use $L^{\infty}$ norm in \cref{m-holder}, the subscript is omitted, i.e., we write $f\in C^{k,\omega}(x_0,t_0)$.
%\end{remark}

Furthermore, we define some other types of continuity.
\begin{definition}\label{d-f2}
Let $U,f$ and $\alpha$ be as in \Cref{d-f1}. We say that $f$ is $C_{s,s'}^{-k,\alpha}$ ($k\geq 1$) at $X_0\in U$ or $f\in C_{s,s'}^{-k,\alpha}(X_0)$ if there exist $K,r_0>0$ such that $Q_{r}(X_0)\subset U$ and
\begin{equation}\label{e.c-1}
\|f\|^*_{L^{s,s'}(Q_r(X_0) )}\leq K r^{-k+\alpha}, ~\forall ~0<r<r_0,
\end{equation}
where we require
\begin{equation}\label{e.c-2}
\frac{n}{s}+\frac{2}{s'}-k+\alpha\geq 0.
\end{equation}
Then define
\begin{equation*}
\|f\|_{C_{s,s'}^{-k,\alpha}(X_0)}=\min \left\{K \big | \cref{e.c-1} ~\mbox{holds with}~K\right\}.
\end{equation*}
If $f\in C_{s,s'}^{-k, \alpha}(X)$ for any $X\in U'\subset U$ with the same $r_0$ and
\begin{equation*}
  \|f\|_{C_{s,s'}^{-k,\alpha}(\bar{U}')}:= \sup_{X\in U'} \|f\|_{C_{s,s'}^{-k,\alpha}(X)}<+\infty,
\end{equation*}
we say that $f\in C_{s,s'}^{-k,\alpha}(\bar{U}')$.

%For the special case, if $n/s+2/s'-k+\alpha=0$, we say
%\begin{equation*}
%\|f\|_{C_{s,s'}^{-k,\alpha}(x_0,t_0)}=\|f\|^*_{L^{s,s'}(U\cap Q_{r_0}(x_0,t_0))}.
%\end{equation*}
\end{definition}

\begin{remark}\label{r-df2.2}
If $s=s'$ in \Cref{d-f1} and \Cref{d-f2}, we may write $f\in C_s^{k,\alpha}(X_0)$ for short.  In this paper, if we use these two definitions for some function, we always assume that $r_0=1$ without loss of generality.
\end{remark}

\begin{remark}\label{r-df1.5}
The requirement \cref{e.c-2} is to ensure that \cref{e.c-1} makes sense. If $n/s+2/s'-k+\alpha<0$, \cref{e.c-1} implies
\begin{equation*}
\|f\|_{L^{s,s'}(Q_r(X_0) )}\leq K r^{n/s+2/s'-k+\alpha}, ~\forall ~0<r<r_0,
\end{equation*}
which is trivial since the left-hand tends to zero and the right-hand tends to infinity as $r\to 0$. By noting \cref{e1.3} and \cref{e1.4}, the following notations make sense:
\begin{equation*}
\pmb{b},\pmb{d}\in C_{2p,2p'}^{-1,\alpha}(0),~~c\in C_{p,p'}^{-2,\alpha}(0),~~f\in C_{q,q'}^{-2,\alpha}(0),~~\pmb{f}\in C_{2}^{-1,\alpha}(0),
\end{equation*}
which are used in our theorems.
\end{remark}

\begin{remark}\label{r-df1.3}
Since $\pmb{b},\pmb{d}\in L^{2p,2p'}(Q_1), c\in L^{p,p'}(Q_1)$, $f\in L^{q,q'}(Q_1)$ and $\pmb{f}\in L^2(Q_1)$ where $p,p',q$ and $q'$ are fixed in this paper, the subscripts are omitted when we state the pointwise smoothness (in \Cref{d-f1} and \Cref{d-f2}) of $\pmb{b},c,\pmb{d},f$ and $\pmb{f}$, i.e., we write $\pmb{b}\in C^{k,\alpha}(0)$ in place of $\pmb{b}\in C^{k,\alpha}_{2p,2p'}(0)$ etc.

In addition, since $\pmb{a}$ is bounded, we can use any $L^{s,s'}$ norm to describe its pointwise smoothness. In this paper, we use the $L^1$ norm. Moreover, we write $\pmb{a}\in C^{k,\alpha}(0)$ to replace $\pmb{a}\in C_1^{k,\alpha}(0)$.
\end{remark}

\begin{remark}\label{r-df1.4}
Throughout this paper, if we say $f\in C_{s,s'}^{k, \alpha}(X_0)$, we always use $P^f$ to denote the corresponding polynomial in \cref{m-holder}. If $k<0$, we set $P^f\equiv 0$.
\end{remark}
~\\

Our main results are stated in the following. Since we consider the pointwise regularity, we focus on the regularity at $0$ throughout this paper. We use the $L^2$ norm in \Cref{m-holder} for the solution $u$ and we write $u\in C^{k,\alpha}$ instead of $u\in C_2^{k,\alpha}$ for simplicity. The first is the $C^{\alpha}$ regularity.
\begin{theorem}[\textbf{$C^{\alpha}$ regularity}]\label{t-Ca}
Let $0<\alpha<1$ and $u$ be a weak solution of \Cref{e.div}. Suppose that
\begin{equation*}
  \begin{aligned}
& \|\pmb{a}-\pmb{a}_{Q_r}\|^*_{L^{1}(Q_r)}\leq \delta,~\forall ~0<r<1,~~\pmb{b}\in L^{2p,2p'}(Q_1), ~~c\in C^{-2,\alpha}(0),~~\pmb{d}\in C^{-1,\alpha}(0),\\
&f \in C^{-2,\alpha}(0),~~\pmb{f} \in C^{-1,\alpha}(0),
  \end{aligned}
\end{equation*}
where $\pmb{a}_{Q_r}=\fint_{Q_r} \pmb{a}$ and $\delta>0$ (small) depends only on $n,\lambda, \Lambda, p,q$ and $\alpha$.

Then $u\in C^{\alpha}(0)$, i.e., there exists a constant $P$ such that
\begin{equation}\label{e.Ca.esti}
  \|u-P\|^*_{L^2(Q_r)}\leq C r^{\alpha}\left(\|u\|^*_{L^{2}(Q_1)}+\|f\|_{C^{-2,\alpha}(0)}+\|\pmb{f}\|_{C^{-1,\alpha}(0)}\right), ~~\forall ~0<r<1
\end{equation}
and
\begin{equation*}
  |P|\leq C\left(\|u\|^*_{L^{2}(Q_1)}+\|f\|_{C^{-2,\alpha}(0)}+\|\pmb{f}\|_{C^{-1,\alpha}(0)}\right),
\end{equation*}
where $C$ depends only on $n,\lambda, \Lambda, p,q$, $\alpha$, $\pmb{b}$, $\|c\|_{C^{-2,\alpha}(0)}$ and $\|\pmb{d}\|_{C^{-1,\alpha}(0)}$.
\end{theorem}

\begin{remark}\label{r-1.4}
 The smallness condition $\|\pmb{a}-\pmb{a}_{Q_r}\|^*_{L^{1}(Q_r)}\leq \delta$ is the so-called small BMO condition, which has been used widely (see \cite{MR2456271, MR2399163, MR2680179, MR2601069, MR2771670}).
\end{remark}

\begin{remark}\label{r-1.5}
We use the $L^1$ norm for $\pmb{a}$, which is adequate for the proof since $\pmb{a}\in L^{\infty}$. This observation is inspired by \cite[P. 418]{MR3620893}.
\end{remark}

Since the pointwise regularity implies the local regularity, we have
\begin{corollary}[\textbf{$C^{\alpha}$ regularity}]\label{co1.0}
Let $0<\alpha<1$ and $u$ be a weak solution of \Cref{e.div}. Suppose that
\begin{equation*}
  \begin{aligned}
& \|\pmb{a}-\pmb{a}_{Q_r(X_0)}\|^*_{L^{1}(Q_r(X_0))}\leq \delta,~\forall ~X_0\in Q_{1/2}, ~\forall~ 0<r<1/2,~~\\
&\pmb{b}\in L^{2p,2p'}(Q_{3/4}), ~~c\in C^{-2,\alpha}(\bar{Q}_{3/4}),~~\pmb{d}\in C^{-1,\alpha}(\bar{Q}_{3/4}),\\
&f \in C^{-2,\alpha}(\bar{Q}_{3/4}),~~\pmb{f} \in C^{-1,\alpha}(\bar{Q}_{3/4}),
  \end{aligned}
\end{equation*}
where $\delta>0$ (small) depends only on $n,\lambda,\Lambda, p,q$ and $\alpha$.

Then $u\in C^{\alpha}(\bar{Q}_{1/2})$ and
\begin{equation*}
  \|u\|_{C^{\alpha}(\bar{Q}_{1/2})}\leq C\left(\|u\|_{L^{2}(Q_{3/4})}+\|f\|_{C^{-2,\alpha}(\bar{Q}_{3/4})}
  +\|\pmb{f}\|_{C^{-1,\alpha}(\bar{Q}_{3/4})}\right),
\end{equation*}
where $C$ depends only on $n,\lambda, \Lambda, p,q$, $\alpha$, $\pmb{b}$, $\|c\|_{C^{-2,\alpha}(\bar{Q}_{3/4})}$ and $\|\pmb{d}\|_{C^{-1,\alpha}(\bar{Q}_{3/4})}$.
\end{corollary}
\begin{remark}\label{re1.4}
In \Cref{d-f1} and \Cref{d-f2}, we only give the definition of $C^{k,\alpha}(\bar{U}')$ for $U'\subset Q_1$. Hence, we use $c\in C^{-2,\alpha}(\bar{Q}_{3/4})$ instead of $c\in C^{-2,\alpha}(\bar{Q}_{1})$. Of course, we can define the later notation and use it in \Cref{co1.0}. Since we confine us to the interior regularity, do not do that.
\end{remark}
~\\

%\begin{remark}\label{r-1.5}
%The smallness condition for $a^{ij}$ is necessary, which is different from the boundary regularity for equations in non-divergence form \cite{ours}. The reason is that $x_n$ is a solution for the boundary situation for equations in non-divergence form.
%\end{remark}

%\begin{theorem}[\textbf{$C^{1,\alpha}$ regularity}]\label{t-C1a}
%Let $0<\alpha<1$ and $u$ be a weak solution of \Cref{e.div}. Suppose that
%\begin{equation*}
%  \begin{aligned}
%&\pmb{a}\in C_2^{\alpha}(0),~~\pmb{b}\in C^{-1,\alpha}(0),
%~~c\in C^{-1,\alpha}(0),~~\pmb{d}\in C^{\alpha}(0),\\
%&f \in C^{-1,\alpha}(0),~~\pmb{f} \in C^{\alpha}(0).
%  \end{aligned}
%\end{equation*}
%
%Then $u\in C_2^{1,\alpha}(0)$, i.e., there exists $P\in \mathcal{P}_1$ such that
%\begin{equation}\label{e.C1a.esti}
%  \|u-P\|^*_{L^2(Q_r)}\leq C r^{1+\alpha}\left(\|u\|^*_{L^{2}(Q_1)}+\|f\|_{C^{-1,\alpha}(0)}
%  +\|\pmb{f}\|_{C^{\alpha}(0)}\right), ~~\forall ~0<r<1,
%\end{equation}
%and
%\begin{equation*}
%  |P(0)|+|DP(0)|\leq C\left(\|u\|^*_{L^{2}(Q_1)}+\|f\|_{C^{-1,\alpha}(0)}
%  +\|\pmb{f}\|_{C^{\alpha}(0)}\right),
%\end{equation*}
%where $C$ depends only on $n,\lambda, \Lambda,p,q,\alpha$,$\|\pmb{a}\|_{C_{2}^{\alpha}(0)}$ $\|\pmb{b}\|_{C^{-1,\alpha}(0)}$, $\|c\|_{C^{-1,\alpha}(0)}$ and $\|\pmb{d}\|_{C^{\alpha}(0)}$.
%\end{theorem}

Now, we state the higher order regularity.

\begin{theorem}\label{t-Cka2}
Let $0<\alpha<1$ and $u$ be a weak solution of \Cref{e.div}. Suppose that for some integers $k\geq 0$ and $l\geq 1$, $u\in C^{k,\alpha}(0)$,
\begin{equation*}
u(0)=|Du(0)|=\cdots=|D^ku(0)|=0,
\end{equation*}
and
\begin{equation*}
  \begin{aligned}
    &\pmb{a}\in C^{l-1,\alpha}(0),~~\pmb{b}\in C^{l-2,\alpha}(0),~~c\in C^{l-3,\alpha}(0),~~\pmb{d}\in C^{l-2,\alpha}(0),\\
    & f\in C^{k+l-2,\alpha}(0),~~\pmb{f}\in C^{k+l-1,\alpha}(0).
  \end{aligned}
\end{equation*}

Then $u\in C^{k+l,\alpha}(0)$, i.e., there exists $P\in \mathcal{P}_{k+l}$ such that
\begin{equation*}
  \begin{aligned}
\|u-P\|^*_{L^2(Q_r)}\leq C r^{k+l+\alpha}\left(\|u\|^*_{L^2(Q_1)}
+\|f\|_{C^{k+l-2,\alpha}(0)}+\|\pmb{f}\|_{C^{k+l-1,\alpha}(0)}\right),~\forall 0<r<1,
  \end{aligned}
\end{equation*}

\begin{equation}\label{e1.12}
\begin{aligned}
\mathbf{\Pi}_{k+l-2} \left(P_{t}-(P^{a^{ij}}P_{i}+P^{d^j}P)_j+P^{b^i} P_{i}+P^cP-P^{f}+P^{f^i}_i\right)=0
\end{aligned}
\end{equation}
and
\begin{equation*}
  |D^{k+1}P(0)|+\cdots+|D^{k+l}P(0)|\leq C\left(\|u\|^*_{L^2(Q_1)}
+\|f\|_{C^{k+l-2,\alpha}(0)}+\|\pmb{f}\|_{C^{k+l-1,\alpha}(0)}\right),
\end{equation*}
where $C$ depends only on $n,\lambda, \Lambda, p, q, \alpha$, $\|\pmb{a}\|_{C^{l-1,\alpha}(0)}$, $\|\pmb{b}\|_{C^{l-2,\alpha}(0)}$, $\|c\|_{C^{l-3,\alpha}(0)}$ and $\|\pmb{d}\|_{C^{l-2,\alpha}(0)}$.
\end{theorem}

\begin{remark}\label{re.kl.1}
This theorem shows that the smooth conditions for the coefficients $\pmb{a},\pmb{b}$ (resp. $c,\pmb{d}$) can be reduced by order $k$ (resp. $k+1$) if $u(0)=\cdots=|D^k u(0)|=0$. This phenomenon was first observed by Han \cite{MR1658901}, who proved the case $l=1$ under the assumption $\pmb{b},c\in L^{\infty}$. We are inspired by \cite{MR4088470}, where a similar assertion for the boundary was proved. Similar boundary regularity were systematically developed in \cite{MR4644419} for the Dirichlet problem and the oblique derivative problem.
\end{remark}

\begin{remark}\label{re1.7}
In fact, the assumption ``$u\in C^{k,\alpha}(0)$'' can be removed. If $u(0)=0$, by \Cref{t-Cka2}, $u\in C^{1,\alpha}(0)$. Thus, $|Du(0)|=0$ makes sense. Once $|Du(0)|=0$ is assumed, $u\in C^{2,\alpha}(0)$. Hence, $|D^2u(0)|=0$ makes sense...
\end{remark}
~\\

As a corollary, we have
\begin{theorem}[\textbf{$C^{k,\alpha}$ regularity}]\label{t-Cka}
Let $k\geq 1$, $0<\alpha<1$ and $u$ be a weak solution of \Cref{e.div}. Suppose that
\begin{equation*}
  \begin{aligned}
&\pmb{a}\in C^{k-1,\alpha}(0),~~\pmb{b}\in C^{k-2,\alpha}(0),
~~c\in C^{k-2, \alpha}(0),~~\pmb{d}\in C^{k-1,\alpha}(0),\\
&f \in C^{k-2,\alpha}(0),~~\pmb{f}\in C^{k-1,\alpha}(0).
  \end{aligned}
\end{equation*}

Then $u\in C^{k,\alpha}(0)$, i.e., there exists $P\in \mathcal{P}_k$ such that
\begin{equation}\label{e.Cka.esti}
  \|u-P\|^*_{L^2(Q_r)}\leq C r^{k+\alpha}\left(\|u\|^*_{L^{2}(Q_1)}+\|f\|_{C^{k-2,\alpha}(0)}
  +\|\pmb{f}\|_{C^{k-1,\alpha}(0)}\right), ~~\forall ~0<r<1,
\end{equation}
\begin{equation*}
  \mathbf{\Pi}_{k-2} \left(P_{t}-(P^{a^{ij}}P_{i}+P^{d^j}P)_j+P^{b^i} P_{i}+P^cP-P^{f}+P^{f^i}_i\right)=0
\end{equation*}
and
\begin{equation*}
  |P(0)|+|DP(0)|+\cdots+|D^k P(0)|\leq C\left(\|u\|^*_{L^{2}(Q_1)}+\|f\|_{C^{k-2,\alpha}(0)}
  +\|\pmb{f}\|_{C^{k-1,\alpha}(0)}\right),
\end{equation*}
where $C$ depends only on $n,\lambda, \Lambda, p, q, \alpha$, $\|\pmb{a}\|_{C^{k-1,\alpha}(0)}$, $\|\pmb{b}\|_{C^{k-2,\alpha}(0)}$, $\|c\|_{C^{k-2,\alpha}(0)}$ and $\|\pmb{d}\|_{C^{k-1,\alpha}(0)}$.
\end{theorem}
\begin{remark}\label{re1.6}
To the best of our knowledge, \Cref{t-Cka} is new. Furthermore, this result is sharp.
\end{remark}
~\\

Similarly, we have the following local regularity.
\begin{corollary}[\textbf{$C^{k,\alpha}$ regularity}]\label{co1.1}
Let $k\geq 1, 0<\alpha<1$ and $u$ be a weak solution of \Cref{e.div}. Suppose that
\begin{equation*}
  \begin{aligned}
& \pmb{a}\in C^{k-1,\alpha}(\bar{Q}_{3/4}),~~\pmb{b}\in C^{k-2,\alpha}(\bar{Q}_{3/4}), ~~c\in C^{k-2,\alpha}(\bar{Q}_{3/4}),~~\pmb{d}\in C^{k-1,\alpha}(\bar{Q}_{3/4}),\\
&f \in C^{k-2,\alpha}(\bar{Q}_{3/4}),~~\pmb{f} \in C^{k-1,\alpha}(\bar{Q}_{3/4}).
  \end{aligned}
\end{equation*}

Then $u\in C^{k,\alpha}(\bar{Q}_{1/2})$ and
\begin{equation*}
  \|u\|_{C^{k,\alpha}(\bar{Q}_{1/2})}\leq C\left(\|u\|_{L^{2}(Q_{3/4})}+\|f\|_{C^{k-2,\alpha}(\bar{Q}_{3/4})}
  +\|\pmb{f}\|_{C^{k-1,\alpha}(\bar{Q}_{3/4})}\right),
\end{equation*}
where $C$ depends only on $n,\lambda, \Lambda, p,q$, $\alpha$, $\|\pmb{b}\|_{C^{k-2,\alpha}(\bar{Q}_{3/4})}$, $\|c\|_{C^{k-2,\alpha}(\bar{Q}_{3/4})}$ and $\|\pmb{d}\|_{C^{k-1,\alpha}(\bar{Q}_{3/4})}$.
\end{corollary}
~\\

As applications, we have the following characterization of the structure of the nodal sets of solutions. For
$k\geq 1$, define
\begin{equation*}
  \mathcal{L}_k(u)=\left\{(x,t)\in \bar Q_{1/2}:~~u(x,t)=\cdots=|D^{k-1}(x,t)|=0,~~D^ku(x,t)\neq 0 \right\}.
\end{equation*}
Then we have
\begin{theorem}\label{th1.3}
Let $l\geq 1$, $0<\alpha<1$ and $u$ be a weak solution of \Cref{e.div}. Suppose that
\begin{equation*}
  \begin{aligned}
    &\pmb{a}\in C^{l-1,\alpha}(\bar{Q}_{3/4}),~~\pmb{b}\in C^{l-2,\alpha}(\bar{Q}_{3/4}),~~c\in C^{l-3,\alpha}(\bar{Q}_{3/4}),~~\pmb{d}\in C^{l-2,\alpha}(\bar{Q}_{3/4}).
  \end{aligned}
\end{equation*}
Then
\begin{equation}\label{e1.6}
  \begin{aligned}
\mathcal{L}_1(u)=\bigcup_{j=0}^{n+1} \mathcal{L}_1^j, \quad
\mathcal{L}_k(u)=\bigcup_{j=0}^{n} \mathcal{L}_k^j~~(k\geq 2),\\
  \end{aligned}
\end{equation}
where each $\mathcal{L}_k^j$ is on a finite union of $j$-dimensional $C^{l,\alpha}$ manifolds.
\end{theorem}
\begin{remark}\label{re1.3}
Note that in above theorem, we have used the parabolic Hausdorff dimension (see \cite[P. 783-784]{MR673830}). We use  the Whitney extension theorem (parabolic version) to prove above result. As far as we know, there is no reference for this result. Hence, we prove it in the Appendix (see \Cref{S.A}).
\end{remark}

\begin{remark}\label{re1.5}
Han \cite{MR1794573} obtained the elliptic version of \Cref{th1.3} (see \Cref{th1.4}). He assumed $\pmb{b},c\in L^{\infty}$ and proved that the manifolds are $C^{1,\beta}$ for \emph{some} $0<\beta<1$. This result corresponds to \Cref{th1.3} (\Cref{th1.4}) with $\pmb{b}\in C^{-1,\alpha},c\in C^{-2,\alpha}$ and we obtain that the manifolds are exactly $C^{1,\alpha}$. Moreover, \Cref{th1.3} shows that if the coefficients have higher regularity, the manifolds are smoother.
\end{remark}
~\\

This paper is organized as follows. In \Cref{sec:1}, we first provide the compactness results for parabolic equations and then present the proof of $C^{\alpha}$ regularity. We use induction to prove \Cref{t-Cka2} and the case of $k=0, l=1$ is proved in \Cref{sec:2}. In \Cref{sec:3}, we give the proof of \Cref{t-Cka2} for $k\geq 0,l\geq 1$ and \Cref{t-Cka} can be obtained immediately by a normalization. The analogous theorems for elliptic equations are provided in \Cref{sec:4}.

\begin{notation}\label{no1.1}
\begin{enumerate}~~\\
\item $\{e_i\}^{n}_{i=1}$: the standard basis of $\mathbb{R}^n$, i.e., $e_i=(0,...0,\underset{i^{th}}{1},0,...0)$.
\item $x=(x_1,x_2,...,x_n) \in \mathbb{R}^{n}$ and $X=(x,t)\in \mathbb{R}^{n+1}$.
\item $|x|:=\left(\sum_{i=1}^{n} x_i^2\right)^{1/2}$ for $x\in \mathbb{R}^n$  and $|X|=|(x,t)|:=(|x|^2+|t|)^{1/2}$ for $X=(x,t)\in \mathbb{R}^{n+1}$.
\item $B_r(x_0)=B(x_0,r)=\{x\in \mathbb{R}^{n}\big| |x-x_0|<r\}$ and $B_r=B_r(0)$.
\item $Q_r(X_0)=Q(X_0,r)=B_r(x_0)\times (t_0-r^2,t_0]$ for $X_0=(x_0,t_0)$ and $Q_r=Q_r(0)$.
\item $Q'_r(X_0)=Q'(X_0,r)=B_r(x_0)\times (t_0-r^2,t_0+r^2)$ and $Q'_r=Q'_r(0)$.
\item $\bar A $: the closure of $A$; $\mathrm{diam}(A)$: the diameter of $A$, where $ A\subset \mathbb{R}^{n+1}$.
\item We use multi-index notation. Let
\begin{equation*}
\sigma=(\sigma_1,...,\sigma_{n+1})=(\sigma',\sigma_{n+1})=(\chi,\tau)=(\chi_1,...,\chi_n,\tau)\in \mathbb{N}^{n+1}.
\end{equation*}
Define
\begin{equation*}
\begin{aligned}
&|\sigma'|=\sum_{i=1}^{n}\sigma_i, \quad |\sigma|=|\sigma'|+2\sigma_{n+1},\quad\sigma!=\prod_{i=1}^{n+1}(\sigma_i!),\quad
(x,t)^{\sigma}= \prod_{i=1}^{n} x_i^{\sigma_i}\cdot t^{\sigma_{n+1}},\quad\\
&D^{\sigma}\varphi = \frac{\partial^{|\sigma'|+\sigma_{n+1}} \varphi }{\partial x_1^{\sigma_1}\cdots \partial x_n^{\sigma_n}\partial t^{\sigma_{n+1}}}, \quad D^k \varphi=\left\{D^{\sigma}\varphi:|\sigma|=k\right\}, \quad |D^k\varphi |
    =\left(\sum_{|\sigma|=k}|D^{\sigma}\varphi|^2\right)^{1/2}.
\end{aligned}
\end{equation*}
\item $\varphi_t=D_t \varphi$, $\varphi _i=D_i \varphi=\partial \varphi/\partial x _{i}$ $(1\leq i\leq n)$, $\varphi _{ij}=D_{ij}\varphi =\partial ^{2}\varphi/\partial x_{i}\partial x_{j}$ $(1\leq i,j\leq n)$.
\item $\mathcal{P}_k (k\geq 0):$ the set of parabolic polynomials of degree less than or equal to $k$. That is, any $P\in \mathcal{P}_k$ can be written as
\begin{equation*}
P(x,t)=\sum_{|\sigma|\leq k}\frac{a_{\sigma}}{\sigma!}(x,t)^{\sigma},
\end{equation*}
where $a_{\sigma}$ are constants. Define
\begin{equation*}
\|P\|= \sum_{|\sigma|\leq k}|a_{\sigma}|.
\end{equation*}
\item $\mathcal{HP}_k (k\geq 0):$ the set of homogeneous polynomials of degree $k$. That is, any $P\in \mathcal{HP}_k$ can be written as
\begin{equation*}
P(x,t)=\sum_{|\sigma|= k}\frac{a_{\sigma}}{\sigma!}(x,t)^{\sigma}.
\end{equation*}
\item $\mathbf{\Pi}_k:$ The projection from $\mathcal{P}_l$ to $\mathcal{P}_k$ for $l\geq k$. That is, if $P\in \mathcal{P}_l$ is written as
\begin{equation*}
P(x,t)=\sum_{|\sigma|\leq l}\frac{a_{\sigma}}{\sigma!}(x,t)^{\sigma},
\end{equation*}
then
\begin{equation*}
\mathbf{\Pi}_kP(x,t)=\sum_{|\sigma|\leq k}\frac{a_{\sigma}}{\sigma!}(x,t)^{\sigma}.
\end{equation*}
If $k<0$, set $\mathbf{\Pi}_kP\equiv 0$.
\end{enumerate}
\end{notation}

\section{$C^{\alpha}$ regularity}\label{sec:1}
The compactness method is applied in this paper. Before proving the regularity, we provide some compactness results in $L^2(U)$ for the weak solutions. First, we recall the energy inequality for the parabolic equations (see \cite[P.143-144]{MR0241822})).
\begin{lemma}\label{th2.1}
Let $U=\Omega\times (T_1,T_2]$ and $u\in L^2(T_1,T_2;H^1(\Omega))\cap L^{\infty}(T_1,T_2;L^2(\Omega))$ be a weak solution of \cref{e.div} in $U$. Then for any $T_1<t'<T_2$ and $\Omega'\subset\subset \Omega$, we have
\begin{equation}\label{e2.1.1}
  \begin{aligned}
\|u\|_{L^2(t',T_2;H^1(\Omega'))}+\|u\|_{L^{\infty}(t',T_2;L^2(\Omega'))}
\leq C\left(\|u\|_{L^2(U)}+\|f\|_{L^{q,q'}(U)}+\|\pmb{f}\|_{L^2(U)}\right),
  \end{aligned}
\end{equation}
where $C$ depends only on $n,\lambda,\Lambda,t',T_1,T_2,\Omega',\Omega$, $\|\pmb{b}\|_{L^{2p,2p'}(U)}$, $\|c\|_{L^{p,p'}(U)}$ and $\|\pmb{d}\|_{L^{2p,2p'}(U)}$.
\end{lemma}
\begin{remark}\label{re.2.1}
If $f=g+h$ in \cref{e.div} and
\begin{equation*}
g\in L^{q,q'},~~h\in L^{p,p'},
\end{equation*}
then the energy inequality holds as well since
\begin{equation*}
\frac{n}{p}+\frac{2}{p'}= 2<\frac{n}{2}+2=\frac{n}{q}+\frac{2}{q'}.
\end{equation*}
\end{remark}
~\\

The next is the estimate for the difference quotient of $u$ in $t$ (see \cite[P. 160, Theorem 4.2 ]{MR0241822}).
\begin{lemma}\label{th2.2}
Let $u\in L^2(T_1,T_2;H^1(\Omega))\cap L^{\infty}(T_1,T_2;L^2(\Omega))$ be a weak solution of \cref{e.div} in $U$. Then for any $T_1<t_1<t_2\leq T_2$ and $\Omega'\subset\subset \Omega$, we have
\begin{equation}\label{e2.2.1}
  \begin{aligned}
\|u(x,t+h)-u(x,t)\|_{L^2(Q_{t_1,t_2,\Omega'})}
\leq Ch^{1/2}\left(\|u\|_{L^2(U)}+\|f\|_{L^{q,q'}(U)}+\|\pmb{f}\|_{L^2(U)}\right),
~\forall ~|h|\leq \frac{t_1-T_1}{2},
  \end{aligned}
\end{equation}
where $C$ depends only on $n,\lambda,\Lambda,t_1,T_1,T_2,\Omega',\Omega$, $\|\pmb{b}\|_{L^{2p,2p'}(U)}$, $\|c\|_{L^{p,p'}(U)}$ and $\|\pmb{d}\|_{L^{2p,2p'}(U)}$.
\end{lemma}

\begin{remark}\label{r-2.1}
In \cref{e2.2.1}, we have taken the zero extension of $u$ in $\Omega\times (T_2,+\infty)$.
\end{remark}
~\\

In the following, we prove the interior pointwise regularity. When we deal with the scaling argument, we will encounter equations in the following form:
\begin{equation}\label{e.div-2}
u_t-(a^{ij}u_i+d^ju)_j+b^iu_i+cu=f+f^c-f^i_i ~~\mbox{in}~~Q_1,
\end{equation}
where $f^c\in L^{p,p'}(Q_1)$. Note that \Cref{th2.1} and \Cref{th2.2} hold as well for \cref{e.div-2} (see \Cref{re.2.1}).

The following is the key step towards interior $C^{\alpha}$ regularity.
\begin{lemma}\label{l-Ca-key}
For any $0<\alpha<1$, there exists $\delta>0$ depending only on $n,\lambda,\Lambda,p,q$ and $\alpha$ such that if $u$ is a weak solution of \cref{e.div-2} with
\begin{equation}\label{e.small}
  \begin{aligned}
    &\|u\|^*_{L^{2}(Q_1)}\leq 1,~~\|f\|^*_{L^{q,q'}(Q_1)}\leq\delta,~~\|f^c\|^*_{L^{p,p'}(Q_1)}\leq\delta,~~
    \|\pmb{f}\|^*_{L^{2}(Q_1)}\leq\delta,\\
    &\|\pmb{a} -\pmb{a}_{Q_1}\|^*_{L^{1}(Q_1)}\leq \delta,~~
    \|\pmb{b}\|^*_{L^{2p,2p'}(Q_1)}\leq \delta,~~
    \|c\|^*_{L^{p,p'}(Q_1)}\leq \delta,~~
    \|\pmb{d}\|^*_{L^{2p,2p'}(Q_1)}\leq \delta.
  \end{aligned}
\end{equation}
Then there exists a constant $\bar P$ such that
\begin{equation}\label{e-lCa-udis}
  \|u-\bar P\|^*_{L^2(Q_{\eta})}\leq \eta^{\alpha}
\end{equation}
and
\begin{equation}\label{e-lca-p}
|\bar P|\leq \bar C,
\end{equation}
where $\bar C$ depends only on $n,\lambda$ and $\Lambda$, and $0<\eta<1/2$ depends also on $\alpha$.
\end{lemma}

\proof We prove the lemma by contradiction. Suppose that the lemma is false. Then there exist $0<\alpha<1$ and sequences of $u_m,\pmb{a}_{m},\pmb{b}_m,c_m,\pmb{d}_m,f_m,f^c_m,\pmb{f}_m$ such that
\begin{equation*}
  u_{m,t}-(a^{ij}_{m}u_{m,i}+d_m^ju_m)_j+b_m^iu_{m,i}+c_mu_m=f_m+f^c_m-f^i_{m,i}~~\mbox{in}~~Q_1
\end{equation*}
with
\begin{equation*}
  \begin{aligned}
    &\|u_m\|^*_{L^{2}(Q_1)}\leq 1,~~\|f_m\|^*_{L^{q,q'}(Q_1)},
    \|f^c_m\|^*_{L^{p,p'}(Q_1)},\|\pmb{f}_m\|^*_{L^{2}(Q_1)}\leq 1/m, \\
    &\|\pmb{a}_m -\pmb{a}_{Q_1}\|^*_{L^{1}(Q_1)},
    \|\pmb{b}_m\|^*_{L^{2p,2p'}(Q_1)}, \|c_m\|^*_{L^{p,p'}(Q_1)},
    \|\pmb{d}_m\|^*_{L^{2p,2p'}(Q_1)}\leq 1/m,
  \end{aligned}
\end{equation*}
and for any constant $P$ with $|P|\leq \bar C$,
\begin{equation}\label{e-lCa-1}
  \|u_m-P\|^*_{L^{2}(Q_{\eta})}> \eta^{\alpha},
\end{equation}
where $\bar C$ is to be specified later and $0<\eta<1/2$ is taken small such that
\begin{equation}\label{e-lCa-2}
\bar C\eta^{1-\alpha}<1/2.
\end{equation}

Note that $u_m$ are uniformly bounded in $L^2(Q_1)$. In addition, from \Cref{th2.1} and \Cref{th2.2}, for any $Q' \subset \subset Q_1$ and $\varepsilon>0$, there exists $h_0>0$ such that for any $0<|h|<h_0$,
\begin{equation*}
  \|\Delta^{h}_i u_m\|_{L^{2}(Q')}\leq \varepsilon, ~\forall ~1\leq i \leq n~~\mbox{and}~~
  \|\Delta^{h}_t u_m\|_{L^{2}(Q')}\leq \varepsilon,
\end{equation*}
where
\begin{equation*}
  \Delta^{h}_i u_m (x,t)=u_m(x+he_i,t)-u_m(x,t)
  ~~\mbox{and}~~\Delta^{h}_t u_m (x,t)=u_m(x,t+h)-u_m(x,t).
\end{equation*}
Hence, by the Fr\'{e}chet–Kolmogorov theorem (or Kolmogorov–Riesz theorem), $\left\{u_m\right\}$ is pre-compact in $L^2(Q')$ (see \cite[Theorem 5, Corollary 8]{MR2734454} or
\cite[Theorem 1]{MR1930200}). Then there exists a subsequence (denoted by $u_m$ again) such that $u_m$ converges in $L^2$ to some function $\bar u$ on compact subsets of $Q_1$. Then, it can be checked easily (note that $\pmb{a}\to \pmb{a}_{Q_1}$ in $L^1$ implies $\pmb{a}\to \pmb{a}_{Q_1}$ in $L^p$ for any $0<p<+\infty$ since $\pmb{a}\in L^{\infty}$) that $\bar u$ is a weak solution of
\begin{equation}\label{e.heat}
  \bar u_t-a_{Q_1}^{ij} \bar u_{ij}=0~~\mbox{in}~~Q_{1}.
\end{equation}

Since $a_{Q_1}^{ij}$ is a constant matrix, $\bar u$ is smooth. By the $C^1$ estimate for \Cref{e.heat}, there exists a constant $\bar{P}$ such that
\begin{equation*}
  \|\bar u-\bar{P}\|^*_{L^2(Q_r)}\leq  \bar C\|\bar u\|^*_{L^2(Q_{3/4})}r=\lim_{m\to \infty} \bar C\|u_m\|^*_{L^2(Q_{3/4})}r\leq \bar Cr, ~~\forall ~0<r<1/2
\end{equation*}
and
\begin{equation*}
  |\bar{P}|\leq \bar C,
\end{equation*}
where $\bar{C}$ depends only on $n,\lambda$ and $\Lambda$.

By noting \cref{e-lCa-2}, we have
\begin{equation}\label{e-lCa-3}
  \|u-\bar{P}\|^*_{L^{2}(Q_{\eta})}\leq \eta^{\alpha}/2.
\end{equation}
However, from \cref{e-lCa-1},
\begin{equation*}
  \|u_m-\bar{P}\|^*_{L^{2}(Q_{\eta})}> \eta^{\alpha}.
\end{equation*}
Let $m\rightarrow \infty$, we have
\begin{equation*}
    \|u-\bar{P}\|^*_{L^{2}(Q_{\eta})}\geq \eta^{\alpha},
\end{equation*}
which contradicts with \cref{e-lCa-3}.  ~\qed~\\

\begin{remark}\label{re-Ca-key1}
This lemma shows that small prescribed data imply a good approximation to $u$ by solutions of the heat equation. In fact, we can obtain much more: $\forall 0<\alpha<1, k\geq 0$, $\exists \delta$ depending only on $n,k, \lambda, \Lambda, p, q$ and $\alpha$, if \Cref{e.small} holds, then there exists a constant $\bar{P}$ such that
\begin{equation}\label{e2.2}
  \|u-\bar P\|^*_{L^2(Q_{\eta})}\leq \eta^{k+\alpha}
\end{equation}
and
\begin{equation*}
|\bar P|\leq \bar C,
\end{equation*}
where $0<\eta<1/2$ and $\bar C$ depend only on $n,k,\lambda$ and $\Lambda$.

However, we can not obtain the $C^{k,\alpha}$ regularity under \Cref{e.small}. Indeed, \cref{e2.2} is only the approximation to $u$ in one scale (i.e., $Q_{\eta}$). To obtain $C^{k,\alpha}$ regularity, \cref{e2.2} must hold in every scale (i.e., $Q_{\eta^m}$ for every $m\geq 1$). The assumptions on the coefficients to guarantee \cref{e2.2} in every scale are exactly the assumptions for the $C^{k,\alpha}$ regularity (see \Cref{t-Ca-scaling}). Since the desired estimate in $Q_{\eta^m}$ is just a scaling version of \Cref{e2.2}, the assumptions on the coefficients are usually minimal. This demonstrates the benefit of studying the pointwise regularity by the method of compactness and perturbation.
\end{remark}

\begin{remark}\label{re-Ca-key2}
From the proof, we know that the desired $\bar{P}$ in \cref{e-lCa-udis} is exactly derived from the heat equation \Cref{e.heat}. So, we call this technique perturbation.
\end{remark}

\begin{theorem}\label{t-Ca-scaling}
Let $0<\alpha<1$ and $u$ be a weak solution of \cref{e.div}. Assume that
\begin{equation}\label{t-Ca-assum}
  \begin{aligned}
    &\|u\|^*_{L^{2}(Q_1)}\leq 1,~~\|f\|_{C^{-2,\alpha}(0)}\leq\delta,
    ~~\|\pmb{f}\|_{C^{-1,\alpha}(0)}\leq \frac{\delta}{2},~~
    \|\pmb{a}-\pmb{a}_{Q_r}\|^*_{L^{1}(Q_r)}\leq \delta,~\forall ~0<r<1,\\
    &\|\pmb{b}\|^*_{L^{2p,2p'}(Q_1)}\leq \delta,
    ~~\|c\|_{C^{-2,\alpha}(0)}\leq \frac{\delta}{C_0},~~
    \|\pmb{d}\|_{C^{-1,\alpha}(0)}\leq \frac{\delta}{2C_0},
  \end{aligned}
\end{equation}
where $\delta$ is the one in \Cref{l-Ca-key} and $C_0$ depends only on $n,\lambda,\Lambda$ and $\alpha$.

Then $u\in C^{\alpha}(0)$, i.e., there exists a constant $P$ such that
\begin{equation*}\label{e-tCa-scal1}
  \|u-P\|^*_{L^{2}(Q_{r})}\leq C r^{\alpha},~\forall ~0<r<1
\end{equation*}
and
\begin{equation*}\label{e-tCa-scal2}
|P|\leq C,
\end{equation*}
where $C$ depends only on $n,\lambda,\Lambda,p,q$ and $\alpha$.
\end{theorem}

\proof To prove that $u$ is $C^{\alpha}$ at $0$, we only need to prove the following. There exist a sequence $P_m$ ($m\geq -1$) such that for all $m\geq 0$,
\begin{equation}\label{e-tCas-u}
\|u-P_m\|^*_{L^{2}(Q_{\eta^{m}})}\leq \eta ^{m\alpha}
\end{equation}
and
\begin{equation}\label{e-tCas-P}
|P_m-P_{m-1}|\leq \bar C\eta^{(m-1)\alpha},
\end{equation}
where $\bar C$ and $\eta$ are as in \Cref{l-Ca-key}. Indeed, \cref{e-tCas-P} implies that $P_m$ is a Cauchy sequence and there exists a constant $P$ such that
\begin{equation*}
  |P_m-P|\leq C\eta^{m\alpha},~\forall ~m\geq 0.
\end{equation*}
For any $0<r<1$, there exists $m\geq 0$ such that $\eta^{m+1}\leq r<\eta^m$. Then
\begin{equation*}
  \begin{aligned}
\|u-P\|^*_{L^{2}(Q_{r})}\leq& C\|u-P\|^*_{L^{2}(Q_{\eta^m})}
  \leq C\|u-P_m\|^*_{L^{2}(Q_{\eta^m})}+C|P-P_m|\\
  \leq& C\eta ^{m\alpha}\leq C\eta ^{(m+1)\alpha}\leq  Cr^{\alpha}.
  \end{aligned}
\end{equation*}

We prove \cref{e-tCas-u} and \cref{e-tCas-P} by induction. For $m=0$, by setting $P_0=P_{-1}\equiv 0$, \crefrange{e-tCas-u}{e-tCas-P} hold clearly. Suppose that the conclusion holds for $m\leq m_0$. We need to prove the conclusion for $m=m_0+1$.

Let $r=\eta ^{m_0}$, $y=x/r$, $s=t/r^2$ and
\begin{equation}\label{e-tCas-v}
  v(y,s)=\frac{u(x,t)-P_{m_0}}{r^{\alpha}}.
\end{equation}
Then $v$ satisfies
\begin{equation}\label{e-tCas-f}
v_s-(\tilde a^{ij}v_i+\tilde d^jv)_j+\tilde b^iv_i+\tilde cv=\tilde f+\tilde{f}^c-\tilde f^i_i~~\mbox{in}~~Q_1,
\end{equation}
where
\begin{equation}\label{e-tCas-new}
  \begin{aligned}
&\tilde{\pmb{a}}(y,s)=\pmb{a}(x,t),~~\tilde {\pmb{b}}(y,s)=r \pmb{b}(x,t),~~\tilde c(y,s)=r^2c(x,t),~~\tilde{\pmb{d}}(y,s)=r\pmb{d}(x,t),\\
&\tilde{f}(y,s)=\frac{f(x,t)}{r^{\alpha-2}},~~\tilde{f}^c(y,s)=-\frac{c(x,t)P_{m_0}}{r^{\alpha-2}},~~
\tilde{\pmb{f}}(y,s)=\frac{\pmb{f}(x,t)-\pmb{d}(x,t)P_{m_0}}{r^{\alpha-1}}.
  \end{aligned}
\end{equation}

In the following, we show that \cref{e-tCas-f} satisfies the assumptions of \Cref{l-Ca-key}. First, by \cref{e-tCas-u} and \cref{e-tCas-v},
\begin{equation*}
\|v\|^*_{L^{2}(Q_1)}= r^{-\alpha}\|u-P_{m_0}\|^*_{L^{2}(Q_r)}\leq 1.
\end{equation*}
Next, by \cref{t-Ca-assum} and \cref{e-tCas-new},
\begin{equation*}
  \begin{aligned}
    &\|\tilde{\pmb{a}}-\tilde{\pmb{a}}_{Q_1}\|^*_{L^{1}(Q_1)}
    =\|\pmb{a}-\pmb{a}_{Q_r}\|^*_{L^{1}(Q_r)}\leq \delta,\\
    &\|\tilde {\pmb{b}}\|^*_{L^{2p,2p'}(Q_1)}=r \|\pmb{b}\|^*_{L^{2p,2p'}(Q_r)}
    \leq \|\pmb{b}\|^*_{L^{2p,2p'}(Q_1)}\leq \delta,\\
    &\|\tilde c\|^*_{L^{p,p'}(Q_1)}=r^2 \|c\|^*_{L^{p,p'}(Q_r)}
    \leq r^{\alpha} \|c\|_{C^{-2,\alpha}(0)} \leq \delta,\\
    &\|\tilde{\pmb{d}}\|^*_{L^{2p,2p'}(Q_1)}=r \|\pmb{d}\|^*_{L^{2p,2p'}(Q_r)}
    \leq r^{\alpha} \|\pmb{d}\|_{C^{-1,\alpha}(0)} \leq  \delta.
  \end{aligned}
\end{equation*}
Finally, by \cref{e-tCas-P}, there exists a constant $C_0$ depending only on $n,\lambda,\Lambda$ and $\alpha$ such that $\|P_m\|\leq C_0$ ($\forall~0\leq m\leq m_0$). Thus, by \cref{t-Ca-assum} and \cref{e-tCas-new},
\begin{equation*}
  \begin{aligned}
\|\tilde{f}\|^*_{L^{q,q'}(Q_1)}&= r^{2-\alpha}\|f\|^*_{L^{q,q'}(Q_r)}\leq \|f\|_{C^{-2,\alpha}(0)}
\leq \delta,\\
\|\tilde{f}^c\|^*_{L^{p,p'}(Q_1)}
   &\leq r^{2-\alpha}C_0\|c\|^*_{L^{p,p'}(Q_r)}\leq C_0\|c\|_{C^{-2,\alpha}(0)} \leq \delta,\\
   \|\tilde{\pmb{f}}\|^*_{L^{2}(Q_1)}
   &\leq r^{1-\alpha}\left(\|\pmb{f}\|^*_{L^{2}(Q_r)}
   +C_0\|\pmb{d}\|^*_{L^{2p,2p'}(Q_r)}\right)~(\mbox{note that}~2p,2p'\geq 2)\\
   &\leq \|\pmb{f}\|_{C^{-1,\alpha}(0)}
   +C_0\|\pmb{d}\|_{C^{-1,\alpha}(0)} \leq \delta.
  \end{aligned}
\end{equation*}

Therefore, \cref{e-tCas-f} satisfies the assumptions of \Cref{l-Ca-key} and hence there exists a constant $\bar P$ such that
\begin{equation*}
\begin{aligned}
    \|v-\bar P\|^*_{L^{2}(Q _{\eta})}&\leq \eta ^{\alpha},
\end{aligned}
\end{equation*}
and
\begin{equation*}
|\bar P|\leq \bar C.
\end{equation*}
Let $P_{m_0+1}(x,t)=P_{m_0}(x,t)+r^{\alpha}\bar P$. Then \cref{e-tCas-P} hold for $m_0+1$. Recalling \cref{e-tCas-v}, we have
\begin{equation*}
  \begin{aligned}
\|u-P_{m_0+1}\|^*_{L^{2}(Q_{\eta^{m_0+1}})}
&= \|u-P_{m_0}-r^{\alpha}\bar P\|^*_{L^{2}(Q_{\eta r})}\\
&= \|r^{\alpha}v-r^{\alpha}\bar P\|^*_{L^{2}(Q_{\eta})}\\
&\leq r^{\alpha}\eta^{\alpha}=\eta^{(m_0+1)\alpha}.
  \end{aligned}
\end{equation*}
Hence, \cref{e-tCas-u} holds for $m=m_0+1$. By induction, the proof is completed.\qed~\\

\begin{remark}\label{re.Ca.scal.1}
The requirements for $\pmb{a},\pmb{b}$ are the same as \Cref{l-Ca-key} but the requirements for $c,\pmb{d}$ are higher than that in \Cref{l-Ca-key}. Since $P_{m_0}$ is a constant, $cP_{m_0}$ and $\pmb{d}P_{m_0}$ will appear in the right hand term during the scaling.
\end{remark}
~\\

Now, we give the~\\
\noindent\textbf{Proof of \Cref{t-Ca}.} In fact, the proof is mere a normalization procedure. First, let $y=x/\rho$, $s=t/\rho^2$ and $u_1(y,s)=u(x,t)$. Then $u_1$ satisfies
\begin{equation}\label{e-Ca-u1}
u_{1,s}-(\tilde a^{ij}u_{1,i}+\tilde d^ju_1)_j+\tilde b^iu_{1,i}+\tilde cu_1=f_1-f_{1,i}^i~~\mbox{in}~~Q_1,
\end{equation}
where
\begin{equation*}\label{e-u1-1}
  \begin{aligned}
&\tilde{\pmb{a}}=\pmb{a},~~\tilde {\pmb{b}}=\rho \pmb{b},~~\tilde c=\rho^2c,~~\tilde{\pmb{d}}=\rho \pmb{d},~~f_1=\rho^2 f,~~\pmb{f}_1=\rho \pmb{f}.
  \end{aligned}
\end{equation*}
Then, by direct computation,
\begin{equation*}
  \begin{aligned}
    &\|\tilde{\pmb{a}}-\tilde{\pmb{a}}_{Q_{r}}\|^*_{L^{1}(Q_{r})}
    =\|\pmb{a}-\pmb{a}_{Q_{\rho r}}\|^*_{L^{1}(Q_{\rho r})}
    \leq \delta,~\forall ~0<r<1,\\
    &\|\tilde {\pmb{b}}\|^*_{L^{2p,2p'}(Q_1)}=\rho \|\pmb{b}\|^*_{L^{2p,2p'}(Q_{\rho})}
    =\frac{1}{|Q_1|}\|\pmb{b}\|_{L^{2p,2p'}(Q_{\rho})},\\
    &\|\tilde c\|_{C^{-2,\alpha}(0)}\leq \rho^{\alpha} \|c\|_{C^{-2,\alpha}(0)},
    ~~\|\tilde {\pmb{d}}\|_{C^{-1,\alpha}(0)}
    \leq \rho^{\alpha} \| \pmb{d}\|_{C^{-1,\alpha}(0)}
  \end{aligned}
\end{equation*}
and
\begin{equation*}
  \begin{aligned}
    &\|f_1\|_{C^{-2,\alpha}(0)}
    \leq \rho^{\alpha} \|f\|_{C^{-2,\alpha}(0)},
    ~~\|\pmb{f}_1\|_{C^{-1,\alpha}(0)}
    \leq \rho^{\alpha} \|\pmb{f}\|_{C^{-1,\alpha}(0)}.
  \end{aligned}
\end{equation*}

Set
\begin{equation*}
M=\|u_1\|^*_{L^{2}(Q_1)}+\|f\|_{C^{-2,\alpha}(0)}+\|\pmb{f}\|_{C^{-1,\alpha}(0)}
=\|u\|^*_{L^{2}(Q_{\rho})}+\|f\|_{C^{-2,\alpha}(0)}+\|\pmb{f}\|_{C^{-1,\alpha}(0)}.
\end{equation*}
Let $u_2=u_1/M$. Then $u_2$ satisfies
\begin{equation}\label{e-Ca-u2}
  u_{2,s}-(\tilde a^{ij}u_{2,i}+\tilde d^ju_2)_j+\tilde b^iu_{2,i}+\tilde cu_2=f_2-f_{2,i}^i~~\mbox{in}~~Q_1,
\end{equation}
where $f_2=f_1/M$ and $\pmb{f}_2=\pmb{f}_1/M$. Hence,
\begin{equation*}
\begin{aligned}
&\|u_2\|^*_{L^{2}(Q_1)}\leq \frac{1}{M} \|u_1\|^*_{L^{2}(Q_1)}\leq 1,\\
&\|f_2\|_{C^{-2,\alpha}(0)}
    \leq \frac{1}{M}\|f_1\|_{C^{-2,\alpha}(0)}
    \leq \frac{\rho^{\alpha}}{M} \|f\|_{C^{-2,\alpha}(0)}
    \leq \rho^{\alpha},\\
&\|\pmb{f}_2\|_{C^{-1,\alpha}(0)}
    \leq \frac{1}{M}\|\pmb{f}_1\|_{C^{-1,\alpha}(0)}
    \leq \frac{\rho^{\alpha}}{M}\|\pmb{f}\|_{C^{-1,\alpha}(0)}
    \leq \rho^{\alpha}.
\end{aligned}
\end{equation*}

From above arguments, by taking $\rho$ small enough (depending only on $n, \lambda,\Lambda, \alpha, \pmb{b}$, $\|c\|_{C^{-2,\alpha}(0)}$ and $\|\pmb{d}\|_{C^{-1,\alpha}(0)}$), \cref{e-Ca-u2} satisfies the assumptions of \Cref{t-Ca-scaling}. Therefore, $u_2$ and thus $u$ is $C^{\alpha}$ at $0$ and the estimate \cref{e.Ca.esti} holds. \qed~\\

\section{Proof of \Cref{t-Cka2} for $k=0, l=1$}\label{sec:2}
In the following two sections, we prove \Cref{t-Cka2} by induction. In this section, we first prove \Cref{t-Cka2} for $k=0, l=1$. The next section is devoted to general $k, l$.

\begin{lemma}\label{l-C1a-key}
For any $0<\alpha<1$, there exists $\delta>0$ depending only on $n,\lambda,\Lambda,p,q$ and $\alpha$ such that if $u\in C^{\alpha}(0)$ is a weak solution of \cref{e.div-2}
with
\begin{equation}\label{e.C1a.1}
  \begin{aligned}
    &\|u\|^*_{L^{2}(Q_1)}\leq 1,~~u(0)=0,
    ~~\|f\|_{C^{-2,\alpha}(0)}\leq\delta,
    ~~\|f^c\|_{C^{-2,\alpha}(0)}\leq\delta,~~
    \|\pmb{f}\|_{C^{-1,\alpha}(0)}\leq\delta,\\
    &\|\pmb{a}-\pmb{a}_{Q_r}\|^*_{L^1(Q_r)}\leq \delta,~~\forall ~0<r<1,~~
    \|\pmb{b}\|^*_{L^{2p,2p'}(Q_1)}\leq \delta,~~\\
    &\|c\|_{C^{-2,\alpha}(0)}\leq \delta,~~
    \|\pmb{d}\|_{C^{-1,\alpha}(0)}\leq \delta,
  \end{aligned}
\end{equation}
then there exists $\bar P\in \mathcal{HP}_1$ such that
\begin{equation}\label{e-lC1a-udis}
  \|u-\bar P\|^*_{L^2(Q_{\eta})}\leq \eta^{1+\alpha}
\end{equation}
and
\begin{equation}\label{e-lc1a-p}
|D \bar P(0)|\leq \bar C,
\end{equation}
where $\bar C$ depends only on $n,\lambda$ and $\Lambda$, and $0<\eta<1/2$ depends also on $\alpha$.
\end{lemma}

\begin{remark}\label{re.C1a.1}
We omit the proof since it is very similar to that of \Cref{l-Ca-key}. The main difference is $\bar P\in \mathcal{HP}_1$ here. As in the proof of \Cref{l-Ca-key}, we have a sequence of $u_m$ and some $u$ with $u_m\to \bar u$ in $L^2$. By the assumption \cref{e.C1a.1}, $u_m$ posses the $C^{\alpha}$ regularity at $0$. Note that $u_m(0)=0$ and hence
\begin{equation*}
  \|u_m\|^*_{L^2(Q_r)}\leq Cr^{\alpha},~~~~\forall~~0<r<1.
\end{equation*}
By taking the limit,
\begin{equation*}
  \|\bar u\|^*_{L^2(Q_r)}\leq Cr^{\alpha},~~~~\forall~~0<r<1.
\end{equation*}
That is, $\bar u(0)=0$. Since $\bar u$ is a solution of \Cref{e.heat}, by the interior $C^2$ estimate and noting $ (0)=0$, there exists $\bar{P}\in \mathcal{HP}_1$ such that
\begin{equation*}
  \|\bar u-\bar{P}\|^*_{L^2(Q_r)}\leq C r^2, ~~\forall ~0<r<1
\end{equation*}
and
\begin{equation*}
  |D \bar{P}(0)|\leq \bar C.
\end{equation*}
The rest proof is almost the same as that of \Cref{l-Ca-key} and we omit it.
\end{remark}
~\\

%\begin{remark}\label{re.C1a.key1}
%The small conditions for the prescribed data are the one in $C^{\alpha}$ regularity. Similarly, the conditions for the prescribed data in the key step of $C^{k,\alpha}$ regularity should be the one in $C^{k-1,\alpha}$ regularity. The reason is that we want to find a $1$-form polynomial to approach $u$.
%\end{remark}

%\begin{remark}\label{re.C1a.key2}
%The condition $u(0)=0$ is needed for \Cref{t-C1a-scaling}. For the higher order regularity ($C^{k,\alpha}$), $u(0)=\cdots=D^{k-1}u(0)=0$ is needed.
%\end{remark}

\begin{theorem}\label{t-C1a-scaling}
Let $0<\alpha<1$ and $u$ be a weak solution of \cref{e.div}. Assume that
\begin{equation}\label{t-C1a-assum}
  \begin{aligned}
    &\|u\|^*_{L^{2}(Q_1)}\leq 1,~~u(0)=0,
    ~~\|f\|_{C^{-1,\alpha}(0)}\leq \delta,~~
    \|\pmb{f}\|_{C^{\alpha}(0)}\leq \frac{\delta}{3},\\
    &\|\pmb{a}-\pmb{a}_{Q_r}\|^*_{L^1(Q_r)}\leq \frac{\delta}{6\sqrt{n}\Lambda C_0}r^{\alpha},~~\forall~~0<r<1,~~
    \|\pmb{b}\|_{C^{-1,\alpha}(0)}\leq \frac{\delta}{2C_0},\\
    &\|c\|_{C^{-2,\alpha}(0)}\leq \frac{\delta}{2C_0},~~
    \|\pmb{d}\|_{C^{-1,\alpha}(0)}\leq \frac{\delta}{3C_0},
  \end{aligned}
\end{equation}
where $\delta$ is the one in \Cref{l-C1a-key} and $C_0$ depends only on $n,\lambda,\Lambda$ and $\alpha$.

Then $u\in C^{1, \alpha}(0)$, i.e., there exists $P\in \mathcal{HP}_1$ such that
\begin{equation*}\label{e-tC1a-scal1}
  \|u-P\|^*_{L^{2}(Q_{r})}\leq C r^{1+\alpha},~\forall ~0<r<1
\end{equation*}
and
\begin{equation*}\label{e-tC1a-scal2}
|DP(0)|\leq C,
\end{equation*}
where $C$ depends only on $n,\lambda,\Lambda$ and $\alpha$.
\end{theorem}

\proof The proof is similar to that of \Cref{t-Ca-scaling} and we only need to prove the following. There exist a sequence of $P_m\in \mathcal{HP}_1$ ($m\geq -1$) such that for all $m\geq 0$,
\begin{equation}\label{e-tC1as-u}
\|u-P_m\|^*_{L^{2}(Q_{\eta^{m}})}\leq \eta ^{m(1+\alpha)}
\end{equation}
and
\begin{equation}\label{e-tC1as-P}
|DP_m(0)-DP_{m-1}(0)|\leq C\eta^{(m-1)\alpha},
\end{equation}
where $\eta$ and $C$ are as in \Cref{l-C1a-key}.

We prove the above by induction. For $m=0$, by setting $P_0=P_{-1}\equiv 0$, \crefrange{e-tC1as-u}{e-tC1as-P} hold clearly. Suppose that the conclusion holds for $m\leq m_0$. We need to prove the conclusion for $m=m_0+1$.

Let $r=\eta ^{m_0}$, $y=x/r$, $s=t/r^2$ and
\begin{equation}\label{e-tC1as-v}
  v(y,s)=\frac{u(x,t)-P_{m_0}(x,t)}{r^{1+\alpha}}.
\end{equation}
Then similar to the previous, $v$ satisfies
\begin{equation}\label{e-tC1as-f}
v_s-(\tilde a^{ij}v_i+\tilde d^jv)_j+\tilde b^iv_i+\tilde cv=\tilde f+\tilde{f}^c-\tilde f^i_i~~\mbox{in}~~Q_1.
\end{equation}
where
\begin{equation}\label{e-tC1as-new}
  \begin{aligned}
&\tilde{\pmb{a}}=\pmb{a},~~\tilde {\pmb{b}}=r \pmb{b},~~\tilde c=r^2c,~~\tilde{\pmb{d}}=r\pmb{d},\\
&\tilde{f}=\frac{f}{r^{\alpha-1}},
~~\tilde{f}^c=-\frac{b^iP_{m_0,i}+cP_{m_0}}{r^{\alpha-1}},\\
&\tilde{f}^i=\frac{(f^i-f^i(0))-(a^{ij}-a^{ij}_{Q_r})P_{m_0,j}-d^{i}P_{m_0}}{r^{\alpha}}.
  \end{aligned}
\end{equation}

As before, we can show that \cref{e-tC1as-f} satisfies the assumptions of \Cref{l-C1a-key}. Indeed, for any $0<\rho<1$,
\begin{equation*}
  \begin{aligned}
    &\|v\|^*_{L^{2}(Q_1)}= r^{-1-\alpha}\|u-P_{m_0}\|^*_{L^{2}(Q_r)}
    \leq 1,\\
    &\|\tilde{\pmb{a}}-\tilde{\pmb{a}}_{Q_{\rho}}\|^*_{L^{1}(Q_{\rho})}
    =\|\pmb{a}-\pmb{a}_{Q_{\rho r}}\|^*_{L^{1}(Q_{\rho r})}\leq \delta,\\
    &\|\tilde {\pmb{b}}\|^*_{L^{2p,2p'}(Q_1)}=r \|\pmb{b}\|^*_{L^{2p,2p'}(Q_{\rho})}
    \leq \delta,\\
    &\|\tilde c\|^*_{L^{p,p'}(Q_{\rho})}=r^{2} \|c\|^*_{L^{p,p'}(Q_{\rho r})}
    \leq r^{\alpha} \|c\|_{C^{-2,\alpha}(0)} \rho^{-2+\alpha}
    \leq \delta \rho^{-2+\alpha},\\
    &\|\tilde {\pmb{d}}\|^*_{L^{2p,2p'}(Q_{\rho})}=r \|\pmb{d}\|^*_{L^{2p,2p'}(Q_{\rho r})}
    \leq r^{\alpha} \|\pmb{d}\|_{C^{-1,\alpha}(0)} \rho^{-1+\alpha}
    \leq \delta \rho^{-1+\alpha}.
  \end{aligned}
\end{equation*}
Thus,
\begin{equation*}
\|\tilde c\|_{C^{-2,\alpha}(0)}\leq \delta,~~
    \|\tilde{\pmb{d}}\|_{C^{-1,\alpha}(0)}\leq \delta.
\end{equation*}
In addition, by \cref{e-tC1as-P}, there exists a constant $C_0$ depending only on $n,\lambda,\Lambda$ and $\alpha$ such that $\|P_m\|\leq C_0$ ($\forall~0\leq m\leq m_0$). Thus, by noting $P_{m_0}\in \mathcal{HP}_1$, for any $0<\rho<1$,
\begin{equation*}
  \begin{aligned}
   \|\tilde{f}\|^*_{L^{q,q'}(Q_{\rho})}
   &= r^{1-\alpha}\|f\|^*_{L^{q,q'}(Q_{\rho r})}
   \leq \|f\|_{C^{-1,\alpha}(0)}\rho^{-1+\alpha}
   \leq \delta \rho^{-2+\alpha},\\
   \|\tilde{f}^c\|^*_{L^{p,p'}(Q_{\rho})}
   &\leq r^{1-\alpha}\left(C_0\|\pmb{b}\|^*_{L^{2p,2p'}(Q_{\rho r})}
   +C_0r\|c\|^*_{L^{p,p'}(Q_{\rho r})}\right)\\
   &\leq C_0\|\pmb{b}\|_{C^{-1,\alpha}(0)}\rho^{-1+\alpha}
    +C_0\|c\|_{C^{-2,\alpha}(0)} \rho^{-2+\alpha}\leq \delta \rho^{-2+\alpha},\\
   \|\tilde{\pmb{f}}\|^*_{L^{2}(Q_{\rho})}
   &\leq r^{-\alpha}\left(\|\pmb{f}-\pmb{f}(0)\|^*_{L^{2}(Q_{\rho r})}
   +C_0\|\pmb{a}-\pmb{a}_{Q_{\rho r}}\|^*_{L^{2}(Q_{\rho r})}
   +C_0r\|\pmb{d}\|^*_{L^{2p,2p'}(Q_{\rho r})}\right)\\
   &\leq r^{-\alpha}\left(\|\pmb{f}-\pmb{f}(0)\|^*_{L^{2}(Q_{\rho r})}
   +2\sqrt{n}\Lambda C_0\|\pmb{a}-\pmb{a}_{Q_{\rho r}}\|^*_{L^{1}(Q_{\rho r})}
   +C_0r\|\pmb{d}\|^*_{L^{2p,2p'}(Q_{\rho r})}\right)\\
   &\leq \left(\frac{\delta}{3}+\frac{\delta}{3}\right)\rho^{\alpha}
   +C_0\|\pmb{d}\|_{C^{-1,\alpha}(0)}\rho^{-1+\alpha}\leq \delta \rho^{-1+\alpha}.
  \end{aligned}
\end{equation*}
Thus,
\begin{equation*}
\|\tilde f\|_{C^{-2,\alpha}(0)}\leq \delta,~~
\|\tilde f^c\|_{C^{-2,\alpha}(0)}\leq \delta,~~
    \|\tilde{\pmb{f}}\|_{C^{-1,\alpha}(0)}\leq \delta.
\end{equation*}

Since \cref{e-tC1as-f} satisfies the assumptions of \Cref{l-C1a-key}, there exists $\bar P \in \mathcal{HP}_1$ such that
\begin{equation*}
\begin{aligned}
    \|v-\bar P\|^*_{L^{2}(Q _{\eta})}&\leq \eta ^{1+\alpha},
\end{aligned}
\end{equation*}
and
\begin{equation*}
|D\bar P(0)|\leq \bar C.
\end{equation*}
Let $P_{m_0+1}(x,t)=P_{m_0}(x,t)+r^{1+\alpha}\bar P(y,s)$. Then \cref{e-tC1as-u} and \cref{e-tC1as-P} hold for $m=m_0+1$. By induction, the proof is completed.\qed~\\

\begin{remark}\label{re.C1a.scal.2}
The condition
\begin{equation}\label{e3.1}
\|\pmb{a}-\pmb{a}_{Q_r}\|^*_{L^1(Q_r)}\leq Cr^{\alpha},~~\forall~~0<r<1,
\end{equation}
is equivalent to $\pmb{a}\in C^{\alpha}(0)$ (see \cite[Chapter III, Theorem 1.2]{MR717034}\cite[Theorem 3.1]{MR1669352}). Indeed, if $\pmb{a}\in C^{\alpha}(0)$,
\begin{equation*}
  \begin{aligned}
\|\pmb{a}-\pmb{a}_{Q_r}\|^*_{L^1(Q_r)}\leq& \|\pmb{a}-\pmb{a}(0)\|^*_{L^1(Q_r)}
+\|\pmb{a}(0)-\pmb{a}_{Q_r}\|^*_{L^1(Q_r)}\\
\leq&2\|\pmb{a}-\pmb{a}(0)\|^*_{L^1(Q_r)}
\leq 2\|\pmb{a}\|_{C^{\alpha}(0)}r^{\alpha},~~\forall~~0<r<1.
  \end{aligned}
\end{equation*}

Conversely, for any $m\geq 0$ ($\rho:=1/2$),
\begin{equation*}
|\pmb{a}_{Q_{\rho^{m+1}}}-\pmb{a}_{Q_{\rho^m}}|\leq |a-\pmb{a}_{Q_{\rho^{m+1}}}|+|a-\pmb{a}_{Q_{\rho^m}}|.
\end{equation*}
By taking the integral average in $Q_{\rho^{m+1}}$ for both sides,
\begin{equation*}
  \begin{aligned}
|\pmb{a}_{Q_{\rho^{m+1}}}-\pmb{a}_{Q_{\rho^m}}|\leq& \|\pmb{a}-\pmb{a}_{Q_{\rho^{m+1}}}\|^*_{L^1(Q_{\rho^{m+1}})}
+2^{n+2}\|\pmb{a}-\pmb{a}_{Q_{\rho^{m}}}\|^*_{L^1(Q_{\rho^{m}})}\\
  \end{aligned}
\end{equation*}
If \cref{e3.1} holds, then
\begin{equation*}
  \begin{aligned}
|\pmb{a}_{Q_{\rho^{m+1}}}-\pmb{a}_{Q_{\rho^m}}|\leq& C\rho^{(m+1)\alpha}+C\rho^{m\alpha}
\leq C\rho^{m\alpha},
  \end{aligned}
\end{equation*}
which implies that $\pmb{a}_{Q_{\rho^m}}$ is a Cauchy sequence. Then there exists a constant matrix $\pmb{a}_0$ such that $\pmb{a}_{Q_{\rho^m}}\to \pmb{a}_0$ and
\begin{equation*}
|\pmb{a}_{Q_{\rho^m}}-\pmb{a}_0|\leq C\rho^{m\alpha}.
\end{equation*}
For any $0<r<1$, there exists $m\geq 0$ such that $\rho^{m+1}\leq r<\rho^{m}$. Then
\begin{equation*}
  \begin{aligned}
\|\pmb{a}-\pmb{a}_{0}\|^*_{L^1(Q_r)}\leq& C\|\pmb{a}-\pmb{a}_{0}\|^*_{L^1(Q_{\rho^{m}})}\\
\leq& C\|\pmb{a}-\pmb{a}_{Q_{\rho^{m}}}\|^*_{L^1(Q_{\rho^{m}})}
+C\|\pmb{a}_{Q_{\rho^{m}}}-\pmb{a}_{0}\|^*_{L^1(Q_{\rho^{m}})}\\
\leq& C\rho^{m\alpha}\leq C\rho^{(m+1)\alpha}\leq Cr^{\alpha},
  \end{aligned}
\end{equation*}
which implies that $\pmb{a}\in C^{\alpha}(0)$.

For simplicity, we use $\pmb{a}\in C^{\alpha}(0)$ in the main theorems \Cref{t-Cka2} and \Cref{t-Cka}.
\end{remark}

\begin{remark}\label{re.C1a.scal.1}
Here, $\pmb{a}\in C^{\alpha}(0)$ and $\pmb{b}\in C^{-1,\alpha}(0)$ are higher than that in \Cref{l-C1a-key}. The reason is that the terms $(a^{ij}-a^{ij}_{Q_r})P_{m_0,j}$ and $P_{m_0,i}b^i$ during the scaling will be treated as right-hand terms. On the contrast, due to $u(0)=0$, the conditions for $c, \pmb{d}$ are the same as in \Cref{l-C1a-key}.
\end{remark}
~\\

Now, we give the~\\
\noindent\textbf{Proof of \Cref{t-Cka2} for $k=0, l=1$.} In fact, the proof is mere a normalization procedure as before. Since the proof is very similar to that of \Cref{t-Ca} and we omit it.\qed~\\

\section{Proof of \Cref{t-Cka2} for $k\geq 0, l\geq 1$}\label{sec:3}
In this section, we first prove \Cref{t-Cka2} for $k\geq 0, l=1$ by induction. For $k=0, l=1$, \Cref{t-Cka2} holds by \Cref{sec:2}. Suppose that \Cref{t-Cka2} holds for $k\leq k_0-1, l=1$. In the following, we prove that \Cref{t-Cka2} holds for $k=k_0, l=1$.

The following is the key step.
\begin{lemma}\label{l-Cka-key}
For any $0<\alpha<1$, there exists $\delta>0$ depending only on $n,k_0,\lambda,\Lambda,p,q$ and $\alpha$ such that if $u\in C^{k_0,\alpha}(0)$ is a weak solution of
\begin{equation}\label{e.div-3}
u_t-(a^{ij}u_i+d^ju)_j+b^iu_i+cu+\hat{P}=f+f^c-f^i_i ~~\mbox{in}~~Q_1,
\end{equation}
where $\hat{P}\in \mathcal{HP}_{k_0-1}$. Assume that
\begin{equation}\label{e.Cka.1}
  \begin{aligned}
    &\|u\|^*_{L^{2}(Q_1)}\leq 1,~~u(0)=|Du(0)|=\cdots=|D^{k_0}u(0)|=0,~~\|\hat{P}\|\leq 1,\\
    &\|f\|_{C^{k_0-2,\alpha}(0)}\leq\delta,~~
    \|f^c\|_{C^{k_0-2,\alpha}(0)}\leq\delta,~~
    \|\pmb{f}\|_{C^{k_0-1,\alpha}(0)}\leq\delta,\\
    &\|\pmb{a}\|_{C^{\alpha}(0)}\leq \delta,~~
    \|\pmb{b}\|_{C^{-1,\alpha}(0)}\leq \delta,~~
    \|c\|_{C^{-2,\alpha}(0)}\leq \delta,~~
    \|\pmb{d}\|_{C^{-1,\alpha}(0)}\leq \delta,
  \end{aligned}
\end{equation}

Then there exists $\bar P\in \mathcal{HP}_{k_0+1}$ such that
\begin{equation}\label{e-lCka-udis}
  \|u-\bar P\|^*_{L^2(Q_{\eta})}\leq \eta^{k_0+1+\alpha},
\end{equation}
\begin{equation*}
  \bar P_t-a^{ij}(0)\bar P_{ij}+\hat{P} \equiv 0
\end{equation*}
and
\begin{equation}\label{e-lCka-p}
|D^{k_0+1} \bar{P}(0)|\leq \bar C,
\end{equation}
where $\bar C$ depend only on $n,k_0,\lambda$ and $\Lambda$,  and $0<\eta<1/2$ depends also on $\alpha$.
\end{lemma}

\begin{remark}\label{re.Cka.1}
We omit the proof since it is very similar to that of \Cref{l-Ca-key} and \Cref{l-C1a-key}. As in the proofs of \Cref{l-Ca-key} and \Cref{l-C1a-key}, we have a sequence of $u_m$ and some $\bar u$ with $u_m\to \bar u$ in $L^2$. By induction and the assumption \cref{e.Cka.1}, $u_m\in C^{k_0,\alpha}(0)$. Note that $D^{k}u_m(0)=0$ for $0\leq k\leq k_0$ and hence
\begin{equation*}
  \|u_m\|^*_{L^2(Q_r)}\leq Cr^{k_0+\alpha},~~~~\forall~~0<r<1.
\end{equation*}
By taking the limit,
\begin{equation*}
  \|\bar u\|^*_{L^2(Q_r)}\leq Cr^{k_0+\alpha},~~~~\forall~~0<r<1.
\end{equation*}
That is, $\bar u(0)=\cdots=|D^{k_0}\bar u(0)|=0$. In addition, we have a sequence of $P_m$ and some $\tilde P$ with $P_m\to \tilde P$. Then $\bar{u}$ is weak solution of
\begin{equation*}
  \bar{u}_t-a^{ij}(0)\bar{u}_{ij}+\tilde P=0~~\mbox{ in}~Q_1.
\end{equation*}

Since $\bar{u}$ is a solution of an equation with constant coefficients, by the interior $C^{k_0+2}$ estimate, there exists $\bar{P}\in \mathcal{HP}_{k_0+1}$ such that
\begin{equation*}
  \|\bar{u}-\bar{P}\|^*_{L^2(Q_r)}\leq C r^{k_0+2}, ~~\forall ~0<r<1,
\end{equation*}
\begin{equation*}
 \bar P_t-a^{ij}(0) \bar P_{ij}+\tilde P\equiv 0
\end{equation*}
and
\begin{equation*}
  |D^{k_0+1} \bar{P}(0)|\leq C.
\end{equation*}
The rest proof is almost the same as before and we omit it.
\end{remark}
~\\

\begin{theorem}\label{t-Cka-scaling}
Let $0<\alpha<1$ and $u\in C^{k_0,\alpha}(0)$ be a weak solution of
\begin{equation}\label{e.div-3.2}
u_t-(a^{ij}u_i+d^ju)_j+b^iu_i+cu+\hat{P}=f-f^i_i ~~\mbox{in}~~Q_1,
\end{equation}
where $P_0\in \mathcal{HP}_{k_0-1}$.
Assume that
\begin{equation}\label{t-Cka-assum}
  \begin{aligned}
    &\|u\|^*_{L^{2}(Q_1)}\leq 1,~~u(0)=|Du(0)|=\cdots=|D^{k_0}u(0)|=0,~~\|\hat{P}\|\leq 1,\\
    &\|f\|^*_{L^{q,q'}(Q_r)}\leq\delta r^{k_0-1+\alpha},~~
    \|\pmb{f}\|^*_{L^2(Q_r)}\leq \frac{\delta}{3}r^{k_0+\alpha},~~\forall~~0<r<1,\\
    &\|\pmb{a}\|_{C^{\alpha}(0)}\leq\frac{\delta}{6\sqrt{n}\Lambda C_0},~~
    \|\pmb{b}\|_{C^{-1,\alpha}(0)}\leq \frac{\delta}{2C_0},\\
    &\|c\|_{C^{-2,\alpha}(0)}\leq \frac{\delta}{2C_0},~~
    \|\pmb{d}\|_{C^{-1, \alpha}(0)}\leq \frac{\delta}{3C_0},
  \end{aligned}
\end{equation}
where $\delta$ is the one in \Cref{l-Cka-key} and $C_0$ depends only on $n,k_0,\lambda,\Lambda$ and $\alpha$.

Then $u\in C^{k_0+1, \alpha}(0)$, i.e., there exists $P\in \mathcal{HP}_{k_0+1}$ such that
\begin{equation*}\label{e-tCka-scal1}
  \|u-P\|^*_{L^{2}(Q_{r})}\leq C r^{k_0+1+\alpha},~\forall ~0<r<1,
\end{equation*}
\begin{equation*}
  P_t-a^{ij}(0)P_{ij}+\hat{P} \equiv 0
\end{equation*}
and
\begin{equation*}\label{e-tCka-scal2}
|D^{k_0+1} P(0)|\leq C,
\end{equation*}
where $C$ depends only on $n,k_0,\lambda,\Lambda$ and $\alpha$.
\end{theorem}

%\begin{remark}\label{re.Cka.2}
%The proof is similar to that of \Cref{t-C1a-scaling} and we only point out the main points here. We prove the theorem by induction.
%\end{remark}
%~\\

\proof As before, we only need to prove the following. There exist a sequence of $P_m\in \mathcal{HP}_{k_0+1}$ ($m\geq 0$) such that for all $m\geq 1$,
\begin{equation}\label{e-tCkas-u}
\|u-P_m\|^*_{L^{2}(Q_{\eta^{m}})}\leq \eta ^{m(k_0+1+\alpha)}
\end{equation}
\begin{equation}\label{e.Cka.2}
  P_{m,t}-a^{ij}(0)P_{m,ij}+\hat{P}\equiv 0
\end{equation}
and
\begin{equation}\label{e-tCkas-P}
|D^{k_0+1} P_m(0)-D^{k_0+1} P_{m-1}(0)|\leq \bar C\eta^{(m-1)\alpha},
\end{equation}
where $\eta$ and $\bar C$ are as in \Cref{l-Cka-key}.

For $m=1$, by \Cref{l-Cka-key} and setting $P_0\equiv 0$, \crefrange{e-tCkas-u}{e-tCkas-P} hold clearly. Suppose that the conclusion holds for $m\leq m_0$. Let $r=\eta ^{m_0}$, $y=x/r$, $s=t/r^2$ and
\begin{equation}\label{e-tCkas-v}
  v(y,s)=\frac{u(x,t)-P_{m_0}(x,t)}{r^{k_0+1+\alpha}}.
\end{equation}
With the aid of \cref{e.Cka.2}, $v$ satisfies
\begin{equation}\label{e-tCkas-f}
v_s-(\tilde a^{ij}v_i+\tilde d^jv)_j+\tilde b^iv_i+\tilde cv=\tilde f+\tilde{f}^c-\tilde f^i_i~~\mbox{in}~~Q_1,
\end{equation}
where
\begin{equation}\label{e-tCkas-new}
  \begin{aligned}
&\tilde {\pmb{a}}=\pmb{a},~~\tilde{ \pmb{b}}=r \pmb{b},~~\tilde c=r^2c,~~\tilde {\pmb{d}}=r\pmb{d},\\
&\tilde{f}=\frac{f}{r^{k_0-1+\alpha}},~~
\tilde{f}^c=-\frac{b^iP_{m_0,i}+cP_{m_0}}{r^{k_0-1+\alpha}},~~
\tilde{f^i}=\frac{f^i-(a^{ij}-a^{ij}(0))P_{m_0,j}-d^iP_{m_0}}{r^{k_0+\alpha}}.
  \end{aligned}
\end{equation}

It can be checked easily that $v,\pmb{a},\pmb{b},c$ and $\pmb{d}$ satisfy the assumptions of \Cref{l-Cka-key}. Next, by \cref{e-tCkas-P}, there exists $C_0$ such that $\|P_m\|\leq C_0$ ($\forall~0\leq m\leq m_0$). Thus, by noting $P_{m_0}\in \mathcal{HP}_{k_0+1}$, for any $0<\rho<1$,
\begin{equation*}
  \begin{aligned}
   \|\tilde{f}\|^*_{L^{q,q'}(Q_{\rho})}
   &= r^{-k_0+1-\alpha}\|f\|^*_{L^{q,q'}(Q_{\rho r})}\leq \delta \rho^{k_0-1+\alpha}
   \leq \delta \rho^{k_0-2+\alpha},\\
   \|\tilde{f}^c\|^*_{L^{p,p'}(Q_{\rho})}
   &\leq r^{-k_0+1-\alpha}\left(C_0r^{k_0}\|\pmb{b}\|^*_{L^{2p,2p'}(Q_{\rho r})}
   +C_0r^{k_0+1}\|c\|^*_{L^{p,p'}(Q_{\rho r})}\right)\\
   &\leq \left(\frac{\delta}{2} \rho^{-1+\alpha}+\frac{\delta}{2} \rho^{-2+\alpha}\right)
   \leq \delta \rho^{k_0-2+\alpha},\\
   \|\tilde{\pmb{f}}\|^*_{L^{2}(Q_{\rho})}
   &\leq r^{-k_0-\alpha}\left(\|\pmb{f}\|^*_{L^{2}(Q_{\rho r})}
   +C_0r^{k_0} \|\pmb{a}-\pmb{a}_{Q_r}\|^*_{L^{2}(Q_{\rho r})}
   +C_0r^{k_0+1} \|\pmb{d}\|^*_{L^{2p,2p'}(Q_{\rho r})}\right)\\
   &\leq (\frac{\delta}{3} \rho^{k_0+\alpha}+\frac{\delta}{3} \rho^{\alpha}+\frac{\delta}{3} \rho^{-1+\alpha})
   \leq \delta \rho^{k_0-1+\alpha}.
  \end{aligned}
\end{equation*}

Since \cref{e-tCkas-f} satisfies the assumptions of \Cref{l-Cka-key} (with $\hat{P}=0$ there), there exists $\bar P\in \mathcal{HP}_{k_0+1}$ such that
\begin{equation*}
\begin{aligned}
    \|v-\bar P\|^*_{L^{2}(Q _{\eta})}&\leq \eta ^{2+\alpha},
\end{aligned}
\end{equation*}
\begin{equation*}
  \bar{P}_t-a^{ij}(0) \bar{P}_{ij}\equiv 0
\end{equation*}
and
\begin{equation*}
|D^{k_0+1}\bar P(0)|\leq C.
\end{equation*}
Let $P_{m_0+1}(x,t)=P_{m_0}(x,t)+r^{k_0+1+\alpha}\bar P(y,s)$. Then \crefrange{e-tCkas-u}{e-tCkas-P} hold for $m=m_0+1$. By induction, the proof is completed.\qed~\\

\begin{remark}\label{re.Cka.key1}
Because of $u(0)=|Du(0)|=\cdots=|D^{k_0}u(0)|=0$, the conditions for $\pmb{a}, \pmb{b}, c$ and $\pmb{d}$ are the same as the ones in $C^{1,\alpha}$ regularity. The reason is that $\pmb{a}, \pmb{b}, c$ and $\pmb{d}$ appear as the right-hand terms in the scaling argument (see \cref{e-tCkas-new}). The $P_{m_0}\in \mathcal{HP}_{k_0+1}$ can provide a quantity $r^{k_0+1}$. Hence, the requirements for these coefficients can be reduced to the ones in $C^{1,\alpha}$ regularity.

On the contrast, we must require $f\in C^{k_0-1,\alpha}(0)$ and $\pmb{f}\in C^{k_0,\alpha}(0)$.
\end{remark}
~\\

Now, we give the~\\
\noindent\textbf{Proof of \Cref{t-Cka2} for $k\geq 0$ and $l=1$.} As before, the proof is mere a normalization procedure. Since we have assumed that \Cref{t-Cka2} holds for $k\leq k_0-1, l=1$, we have $u\in C^{k_0,\alpha}(0)$. By induction, we only need to prove \Cref{t-Cka2} for $k=k_0, l=1$, i.e., $u\in C^{k_0+1,\alpha}(0)$.

Let $\hat{P}=-P^{f}+P^{f^i}_i$. Since $f\in C^{k_0-1,\alpha}(0)$ and $\pmb{f}\in C^{k_0,\alpha}(0)$, $\hat{P}\in \mathcal{P}_{k_0-1}$. In addition, by $u(0)=|Du(0)|\cdots=|D^{k_0}u(0)|=0$ and \cref{e1.12},
\begin{equation*}
\mathbf{\Pi}_{k_0-2} \left(P_{t}-(P^{a^{ij}}P_{i}+P^{d^j}P)_j+P^{b^i} P_{i}+P^cP-P^{f}+P^{f^i}_i\right)=\mathbf{\Pi}_{k_0-2}\hat{P}=0.
\end{equation*}
Thus, $\hat{P}\in \mathcal{HP}_{k_0-1}$.

Add $\hat{P}$ to both sides of \Cref{e.div}. Then $u$ satisfies
\begin{equation*}
u_{t}-(a^{ij}u_{i}+d^ju)_j+b^iu_{i}+cu+\hat{P}=f_1-f_{1,i}^i~~\mbox{in}~~Q_1,
\end{equation*}
where
\begin{equation*}
  \begin{aligned}
    f_1=f-P^{f}~~\mbox{and}~~\pmb {f}_1=\pmb {f}-P^{\pmb{f}}.
  \end{aligned}
\end{equation*}
Hence, for any $0<r<1$,
\begin{equation*}
\|f_1\|^*_{L^{q,q'(Q_r)}}\leq [f]_{C^{k_0-1,\alpha}(0)} r^{k_0-1+\alpha},~~ \|\pmb{f}_1\|^*_{L^{2(Q_r)}}\leq [\pmb{f}]_{C^{k_0,\alpha}(0)} r^{k_0+\alpha}.
\end{equation*}

Next, let $y=x/\rho, s=t/\rho^2$ for $0<\rho<1$ and $u_1(y,s)=u(x,t)$. Then $u_1$ satisfies
\begin{equation*}
u_{1,t}-(\tilde{a}^{ij}u_{1,i}+\tilde{d}^ju_1)_j+\tilde{b}^iu_{1,i}+\tilde{c}u_1+\hat{P}_1=f_2-f_{2,i}^i~~\mbox{in}~~Q_1,
\end{equation*}
where
\begin{equation*}
\begin{aligned}
&\tilde{\pmb{a}}=\pmb{a},~~\tilde {\pmb{b}}=\rho \pmb{b},~~\tilde c=\rho^2c,~~\tilde{\pmb{d}}=\rho \pmb{d},\\
&f_2=\rho^2f_1,~\pmb{f}_2=\rho \pmb{f}_1,~\hat{P}_1=\rho^{2}\hat{P}.
\end{aligned}
\end{equation*}
Let
\begin{equation*}
M=\|u_1\|^*_{L^{2}(Q_1)}+\|f\|_{C^{k_0-1,\alpha}(0)}
+\|\pmb{f}\|_{C^{k_0,\alpha}(0)}
\end{equation*}
and $u_2=u_1/M$. Then $u_2$ satisfies
\begin{equation}\label{e.Cka.5.1}
u_{2,t}-(\tilde{a}^{ij}u_{2,i}+\tilde{d}^ju_2)_j+\tilde{b}^iu_{2,i}+\tilde{c}u_2+\hat{P}_2=f_3-f_{3,i}^i~~\mbox{in}~~Q_1,
\end{equation}
where $f_3=f_2/M$, $\pmb{f}_3=\pmb{f}_2/M$ and $\hat{P}_2=\hat{P}_1/M$.

Therefore, by taking $\rho$ small enough (depending only on $n,k_0,\lambda,\Lambda,p,q,\alpha$, $\|\pmb{a}\|_{C^{\alpha}(0)}$, $\|\pmb{b}\|_{C^{-1,\alpha}(0)}$, $\|c\|_{C^{-2,\alpha}(0)}$ and $\|\pmb{d}\|_{C^{-1,\alpha}(0)}$.), the assumptions \cref{t-Cka-assum} for $u_2$ are satisfied. By \Cref{t-Cka-scaling}, $u_2\in C^{k_0+1,\alpha}(0)$ and thus $u\in C^{k_0+1,\alpha}(0)$.\qed~\\

Now, we give the~\\
\noindent\textbf{Proof of \Cref{t-Cka2} for $k\geq 0$ and $l\geq 1$.} We prove the theorem by induction in $l$.
For $l=1$, the theorem has been proved. Suppose that the theorem holds for $l\leq l_0-1$ and we need to prove it for $l=l_0$. By induction, $u\in C^{k+l_0-1,\alpha}(0)$. That is, there exists $\tilde{P}\in \mathcal{P}_{k+l_0-1}$ such that
\begin{equation*}
\|u-\tilde{P}\|^*_{L^2(Q_r)}\leq C r^{k+l_0-1+\alpha},~\forall 0<r<1
\end{equation*}
and
\begin{equation}\label{e.5.1}
\mathbf{\Pi}_{k+l_0-3} \left(L(\tilde{P})-P^{f}+P^{f^i}_i\right)=0,
\end{equation}
where
\begin{equation*}
L(\tilde{P}):=\tilde{P}_{t}-(P^{a^{ij}}\tilde{P}_{i}+P^{d^j}\tilde{P})_j+P^{b^i} \tilde{P}_{i}+P^c\tilde{P}.
\end{equation*}

Let $v=u-\tilde{P}$. Then
\begin{equation*}
v(0)=\cdots=|D^{k+l_0-1}v(0)|=0
\end{equation*}
and $v$ is a weak solution of
\begin{equation*}
  \begin{aligned}
 &v_{t}-(a^{ij}v_{i}+d^jv)_j+b^iv_{i}+cv\\
&+\tilde{P}_{t}-(a^{ij}\tilde{P}_{i}+d^j\tilde{P})_j+b^i\tilde{P}_{i}+c\tilde{P}=f-f_{i}^i~~\mbox{in}~~Q_1.
  \end{aligned}
\end{equation*}
We rewrite it as
\begin{equation*}
  \begin{aligned}
&v_{t}-(a^{ij}v_{i}+d^jv)_j+b^iv_{i}+cv=\tilde{f}+\tilde{f}^c-\tilde{f}_{i}^i~~\mbox{in}~~Q_1,
  \end{aligned}
\end{equation*}
where
\begin{equation*}
  \begin{aligned}
    &\tilde{f}=f-L(\tilde{P}),\\
    &\tilde{f}^c=-(b^i-P^{b^i})\tilde{P}_{i}-(c-P^c)\tilde{P},\\
    &\tilde{f}^i=f^i+\left(a^{ij}-P^{a^{ij}}\right)\tilde{P}_j+\left(d^i-P^{d^i}\right)\tilde{P}.
  \end{aligned}
\end{equation*}

Since $u(0)=\cdots=|D^ku(0)|=0$, $\mathbf{\Pi}_{k}(\tilde{P})=0$. By combining with the assumptions on the coefficients , for any $0<r<1$,
\begin{equation*}
  \begin{aligned}
&\|\tilde{f}+L(\tilde{P})-P^{f}\|^*_{L^{q,q'}(Q_r)}=\|f-P^f\|^*_{L^{q,q'}(Q_r)}\leq Cr^{k+l_0-2+\alpha},\\
&\|\tilde{f}^c\|^*_{L^{p,p'}(Q_r)}\leq Cr^{k}\|\pmb{b}-P^{\pmb{b}}\|^*_{L^{2p,2p'}(Q_r)}
+Cr^{k+1}\|c-P^{c}\|^*_{L^{p,p'}(Q_r)}\leq Cr^{k+l_0-2+\alpha},\\
&\|\tilde{\pmb{f}}-P^{\pmb{f}}\|^*_{L^2(Q_r)}\leq \|\pmb{f}-P^{\pmb{f}}\|^*_{L^2(Q_r)}+Cr^{k}\|\pmb{a}-P^{\pmb{a}}\|^*_{L^2(Q_r)}+
Cr^{k+1}\|\pmb{d}-P^{\pmb{d}}\|^*_{L^2(Q_r)}\\
&\quad\quad\quad\quad\quad\quad~\leq Cr^{k+l_0-1+\alpha}.\\
  \end{aligned}
\end{equation*}
Hence, by \Cref{t-Cka2} for $k=k+l_0-1$ and $l=1$, there exists $\bar{P}\in \mathcal{HP}_{k+l_0}$ such that
\begin{equation*}
  \|v-\bar{P}\|^*_{L^{2}(Q_{r})}\leq C r^{k+l_0+\alpha},~\forall ~0<r<1,
\end{equation*}
\begin{equation*}
  \begin{aligned}
&\mathbf{\Pi}_{k+l_0-2} \left(L(\bar{P})-P^{\tilde{f}}-P^{\tilde{f}^c}+P^{\tilde{f}^i}_i\right)
=\mathbf{\Pi}_{k+l_0-2} \left(L(\bar{P})+L(\tilde{P})-P^{f}+P^{f^i}_i\right)=0
  \end{aligned}
\end{equation*}
and
\begin{equation*}
|D^{k+l_0} \bar{P}(0)|\leq C,
\end{equation*}

Let $P=\tilde{P}+\bar{P}$. Then
\begin{equation*}
  \|u-P\|^*_{L^{2}(Q_{r})}=\|v-\bar{P}\|^*_{L^{2}(Q_{r})}\leq C r^{k+l_0+\alpha},~\forall ~0<r<1,
\end{equation*}
and
\begin{equation*}
  \begin{aligned}
\mathbf{\Pi}_{k+l_0-2} \left(L(P)-P^{f}+P^{f^i}_i\right)=&\mathbf{\Pi}_{k+l_0-2} \left(L(\bar{P})+L(\tilde{P})-P^{f}+P^{f^i}_i\right)=0.
  \end{aligned}
\end{equation*}
~\qed~\\

Next, we give the~\\
\noindent\textbf{Proof of \Cref{t-Cka}.} Let $v=u-u(0)$. Then $v$ satisfies
\begin{equation*}
  v_{t}-(a^{ij}v_{i}+d^jv)_j+b^iv_{i}+cv=\tilde{f}-\tilde{f}_{i}^i~~\mbox{in}~~Q_1,
\end{equation*}
where $\tilde{f}=f-cu(0)$ and $\tilde{\pmb{f}}=\pmb{f}-\pmb{d}u(0)$. Thus, $v(0)=0$. By $c\in C^{k-2,\alpha}(0)$ and $\pmb{d}\in C^{k-1,\alpha}(0)$, we have
\begin{equation*}
  \tilde{f}\in C^{k-2,\alpha}(0),~~\tilde{\pmb{f}}\in C^{k-1,\alpha}(0).
\end{equation*}
By \Cref{t-Cka2} for $k=0,l=k$, we have $v\in C^{k,\alpha}(0)$ and hence $u\in C^{k,\alpha}(0)$.\qed~\\

Finally, we give the~\\
\noindent\textbf{Proof of \Cref{co1.1}.} Since
\begin{equation*}
  \begin{aligned}
& \pmb{a}\in C^{k-1,\alpha}(\bar{Q}_{3/4}),~~\pmb{b}\in C^{k-2,\alpha}(\bar{Q}_{3/4}), ~~c\in C^{k-2,\alpha}(\bar{Q}_{3/4}),~~\pmb{d}\in C^{k-1,\alpha}(\bar{Q}_{3/4}),\\
&f \in C^{k-2,\alpha}(\bar{Q}_{3/4}),~~\pmb{f} \in C^{k-1,\alpha}(\bar{Q}_{3/4}),
  \end{aligned}
\end{equation*}
by \Cref{a-l2.1}, for any $X_0\in Q_{3/4}$,
\begin{equation*}
  \begin{aligned}
&\pmb{a}\in C^{k-1,\alpha}(X_0),~~\pmb{b}\in C^{k-2,\alpha}(X_0),
~~c\in C^{k-2, \alpha}(X_0),~~\pmb{d}\in C^{k-1,\alpha}(X_0),\\
&f \in C^{k-2,\alpha}(X_0),~~\pmb{f}\in C^{k-1,\alpha}(X_0).
  \end{aligned}
\end{equation*}
Then from \Cref{t-Cka}, $u\in C^{k,\alpha}(X_0)$. That is, there exists $P_{X_0}\in \mathcal{P}_k$ ($k\geq 0$) such that
\begin{equation}\label{a-e3.2-2}
\|u-P_{X_0}\|^*_{L^2(Q(X_0,r))}\leq Cr^{k+\alpha}, ~\forall ~0<r<\frac{1}{4}.
\end{equation}
By \Cref{thA.3}, $u\in C^{k,\alpha}(\bar{Q}_{1/2})$.~\qed~\\

\section{Structure of nodal sets}
In this section, we prove \Cref{th1.3} by the Whitney extension and the implicit function theorem. Our proof is motivated by the structure of singular sets in free boundary problems (see \cite[Theorem 7.9]{MR2962060} and \cite{MR3385162}).
~\\

\noindent\textbf{Proof of \Cref{th1.3}.} Let $E:=\bar{\mathcal{L}}_k\subset \bar{Q}_{1/2}$ be a compact set. We first show that for any $X\in E$, there exists $P_{X}\in \mathcal{P}_{k+l-1}$ such that
\begin{equation}\label{e5.1}
\|u-P_{X}\|^*_{L^{p,p'}(Q(X,r))}\leq Cr^{k+l-1+\alpha}, ~\forall ~0<r<1/2
\end{equation}
and
\begin{equation}\label{e5.2}
|P_X(X)|=\cdots=|D^{k-1}P_X(X)|=0, \quad \|P\|\leq C,
\end{equation}
where $C$ depends only on $n,\lambda, \Lambda, p, q, \alpha$, $\|\pmb{a}\|_{C^{l-1,\alpha}}$, $\|\pmb{b}\|_{C^{l-2,\alpha}}$, $\|c\|_{C^{l-3,\alpha}}$ and $\|\pmb{d}\|_{C^{l-2,\alpha}}$.

Indeed, given $X_0\in E$, there exist a sequence of $X_m\in \mathcal{L}_k$ such that $X_m\to X_0$. By \Cref{t-Cka2}, there exist a sequence of $P_m\in \mathcal{P}_{k+l-1}$ such that \cref{e5.1} and \cref{e5.2} hold for $X_m$. Since $\|P_m\|$ is bounded, up to a subsequence, $P_m\to P_0$ for some $P_0\in \mathcal{P}_{k+l-1}$. By letting $m\to \infty$ in \cref{e5.1} and \cref{e5.2} for $X_m$, we obtain \cref{e5.1} and \cref{e5.2} for $X_0$.

By \Cref{coA.1}, $u$ extends to a $C^{k+l-1,\alpha}$ function in $Q'_1$ (we denote the extension by $u$ again). In the following, we use this extended $u$ instead of the original one. For any $X_0\in \mathcal{L}_k$,
\begin{equation*}
u(X_0)=\cdots=|D^{k-1}u(X_0)|=0, \quad D^{k}u(X_0)\neq 0.
\end{equation*}
Introduce the function $F:Q'_1\to \mathbb{R}^{N}$ (for some integer $N$) by
\begin{equation*}
F(X)=\left\{D^{\sigma}u(X): |\sigma|\leq k-1\right\},~~X\in Q'_1.
\end{equation*}
Since $X_0\in \mathcal{L}_k$,
\begin{equation*}
D_xF(X_0)\neq 0.
\end{equation*}
Otherwise, there must exist $\sigma$ with $|\sigma|=k-2$ such that $D_tD^{\sigma}u(X_0)\neq 0$ (since $X_0\in \mathcal{L}_k$). Then
\begin{equation*}
\mathbf{\Pi}_{k-2} \left(P_{X_0,t}-(P^{a^{ij}}P_{X_0,i}+P^{d^j}P_{X_0})_j+P^{b^i} P_{X_0,i}+P^cP_{X_0}\right)
=P_{X_0,t}\neq 0,
\end{equation*}
which contradicts with \cref{e1.12}.

Note that $D_xF(X_0)$ is a linear operator from $\mathbb{R}^n$ to $\mathbb{R}^{N}$. Let $j$ denote the dimension of the kernel of $D_xF(X_0)$. Since $D_xF(X_0)\neq 0$, $0\leq j\leq n-1$. Then by an rearrangement of the coordinate system, we may assume that
\begin{equation*}
\det\left(\frac{\partial (f^1,...,f^{n-j})}{\partial (x_1,...,x_{n-j})}(X_0)\right)\neq 0,
\end{equation*}
where $f^i$ are components of $F$. By the implicit function theorem, there exists $r>0$ such that
\begin{equation*}
\{(f^1,...,f^{n-j})=0\}\cap Q'_r(X_0)=\left\{(x_1,...,x_{n-j})
=\varphi(x_{n-j+1},...,x_n,t)\right\}\cap Q'_r(X_0)
\end{equation*}
and $\varphi\in C^{l,\alpha}$. Therefore, $\{(f^1,...,f^{n-j})=0\}\cap Q'_r(X_0)$ is a $(j+2)$-dimensional $C^{l,\alpha}$ surface. Note that
\begin{equation*}
\mathcal{L}_k\cap Q'_r(X_0)\subset \{(f^1,...,f^{n-j})=0\}\cap Q'_r(X_0).
\end{equation*}
That is, $\mathcal{L}_k\cap Q'_r(X_0)$ is contained in a $(j+2)$-dimensional $C^{l,\alpha}$ surface.

If there exists $\sigma$ with $|\sigma|=k-2$ such that $D_tD^{\sigma}u(X_0)\neq 0$ in addition, let $\tilde{f}=D^{\sigma}u$. Since $X_0\in \mathcal{L}_k$,
\begin{equation*}
\tilde{f}_t(X_0)\neq 0, \quad \tilde{f}_i (X_0)=0,~\forall ~1\leq i\leq n.
\end{equation*}
Hence,
\begin{equation*}
\det\left(\frac{\partial (f^1,...,f^{n-j},\tilde{f})}{\partial (x_1,...,x_{n-j},t)}(X_0)\right)\neq 0.
\end{equation*}
Therefore, for some $r>0$, $\mathcal{L}_k\cap Q'_r(X_0)$ is contained in a $(j+1)$-dimensional $C^{l,\alpha}$ surface. We remark here that if $j=0$, $X_0$ will be an isolated point. That is,
\begin{equation*}
\mathrm{dim}(\mathcal{L}_k\cap Q'_r(X_0))=\mathrm{dim}(\{X_0\})=0.
\end{equation*}

Finally, we note that if $k\geq 2$ and $j=n-1$, there must exist $\tilde\sigma$ with $|\tilde\sigma|=k-2$ such that $D_tD^{\tilde\sigma}u(X_0)\neq 0$. Otherwise, without loss of generality, we assume
\begin{equation*}
D^{\sigma}u(X_0)=0,~\forall ~\sigma\neq (k,0,0,...,0).
\end{equation*}
Then
\begin{equation*}
\mathbf{\Pi}_{k-2} \left(P_{X_0,t}-(P^{a^{ij}}P_{X_0,i}+P^{d^j}P_{X_0})_j+P^{b^i} P_{X_0,i}+P^cP_{X_0}\right)
=P^{a^{11}}P_{X_0,11}\neq 0,
\end{equation*}
which contradicts with \cref{e1.12}. Let $\tilde{f}=D^{\tilde\sigma}u$ and by an argument as above, $\mathcal{L}_k\cap Q'_r(X_0)$ is contained in a $n$-dimensional $C^{l,\alpha}$ surface for some $r>0$. ~\qed~\\

\section{Elliptic case}\label{sec:4}
In this section, we present the $C^{k,\alpha}$ ($k\geq 0$) regularity for elliptic equation \Cref{e.div-0}. We assume
\begin{equation}\label{e1.2-ell}
 \pmb{b},\pmb{d}\in L^{n}(B_1), c\in L^{\frac{n}{2}}(B_1) ,f\in L^{\frac{2n}{n+2}}(B_1), \pmb{f}\in L^2(B_1).
\end{equation}
Similar to the parabolic case, these are the minimal requirements for the energy inequality (usually called Caccioppoli inequality for elliptic equations). In fact, we can not obtain the energy inequality directly under the assumptions in \cref{e1.2-ell}. We need to assume that the norms are sufficient small, which can be guaranteed by a scaling argument. Hence, \cref{e1.2-ell} is enough.

For \cref{e.div-0}, we have the following interior pointwise regularity.
\begin{theorem}[\textbf{$C^{\alpha}$ regularity}]\label{t-Ca-ell}
Let $0<\alpha<1$ and $u$ be a weak solution of \Cref{e.div-0}. Suppose that
\begin{equation*}
\pmb{b}\in L^{n}(B_1) ~~\mbox{ for}~n\geq 3;~~~~
 |\pmb{b}|^2 \ln \left(1+|\pmb{b}|^2\right)\in L^{1}(B_1)~~\mbox{ for}~n=2
\end{equation*}
and
\begin{equation*}
  \begin{aligned}
&\|\pmb{a}-\pmb{a}_{B_r}\|^*_{L^{1}(B_r)}\leq \delta,~\forall ~0<r<1, ~~c\in C^{-2,\alpha}(0),~~\pmb{d}\in C^{-1,\alpha}(0),\\
&f \in C^{-2,\alpha}(0),~~\pmb{f} \in C^{-1,\alpha}(0),
  \end{aligned}
\end{equation*}
where $\delta>0$ (small) depends only on $n,\lambda, \Lambda$ and $\alpha$.

Then $u\in C^{\alpha}(0)$, i.e., there exists a constant $P$ such that
\begin{equation*}\label{e.Ca.esti-ell}
  \|u-P\|^*_{L^2(B_r)}\leq C r^{\alpha}\left(\|u\|^*_{L^{2}(B_1)}+\|f\|_{C^{-2,\alpha}(0)}+\|\pmb{f}\|_{C^{-1,\alpha}(0)}\right), ~~\forall ~0<r<1
\end{equation*}
and
\begin{equation*}
  |P|\leq C\left(\|u\|^*_{L^{2}(B_1)}+\|f\|_{C^{-2,\alpha}(0)}+\|\pmb{f}\|_{C^{-1,\alpha}(0)}\right),
\end{equation*}
where $C$ depends only on $n,\lambda, \Lambda$, $\alpha$, $\pmb{b}$, $\|c\|_{C^{-2,\alpha}(0)}$ and $\|\pmb{d}\|_{C^{-1,\alpha}(0)}$.
\end{theorem}

\begin{remark}\label{re5.1}
Note that if $n=2$ and $\pmb{b}\in L^2$, we don't have the energy inequality and the $C^{\alpha}$ regularity fails (see counterexamples in \cite{MR3048265}). If we enhance the assumption to $|\pmb{b}|^2 \ln \left(1+|\pmb{b}|^2\right)\in L^{1}(B_1)$, the energy inequality holds (see Lemma 4.2 an (4.2) in \cite{MR3048265}) and therefore we have the $C^{\alpha}$ regularity.
\end{remark}
~\\

\begin{theorem}\label{t-Cka2-ell}
Let $0<\alpha<1$ and $u$ be a weak solution of \Cref{e.div-0}. Suppose that for some integers $k\geq 0$ and $l\geq 1$, $u\in C^{k,\alpha}(0)$,
\begin{equation*}
u(0)=|Du(0)|=\cdots=|D^ku(0)|=0,
\end{equation*}
and
\begin{equation*}
  \begin{aligned}
    &\pmb{a}\in C^{l-1,\alpha}(0),~~\pmb{b}\in C^{l-2,\alpha}(0),~~c\in C^{l-3,\alpha}(0),~~\pmb{d}\in C^{l-2,\alpha}(0),\\
    & f\in C^{k+l-2,\alpha}(0),~~\pmb{f}\in C^{k+l-1,\alpha}(0).
  \end{aligned}
\end{equation*}

Then $u\in C^{k+l,\alpha}(0)$, i.e., there exists $P\in \mathcal{P}_{k+l}$ such that
\begin{equation*}
  \begin{aligned}
\|u-P\|^*_{L^2(B_r)}\leq C r^{k+l+\alpha}\left(\|u\|^*_{L^2(B_1)}
+\|f\|_{C^{k+l-2,\alpha}(0)}+\|\pmb{f}\|_{C^{k+l-1,\alpha}(0)}\right),~\forall 0<r<1,
  \end{aligned}
\end{equation*}

\begin{equation*}\label{e5.12}
\begin{aligned}
\mathbf{\Pi}_{k+l-2} \left(P_{t}-(P^{a^{ij}}P_{i}+P^{d^j}P)_j+P^{b^i} P_{i}+P^cP-P^{f}+P^{f^i}_i\right)=0
\end{aligned}
\end{equation*}
and
\begin{equation*}
  |D^{k+1}P(0)|+\cdots+|D^{k+l}P(0)|\leq C\left(\|u\|^*_{L^2(B_1)}
+\|f\|_{C^{k+l-2,\alpha}(0)}+\|\pmb{f}\|_{C^{k+l-1,\alpha}(0)}\right),
\end{equation*}
where $C$ depends only on $n,\lambda, \Lambda, \alpha$, $\|\pmb{a}\|_{C^{l-1,\alpha}(0)}$, $\|\pmb{b}\|_{C^{l-2,\alpha}(0)}$, $\|c\|_{C^{l-3,\alpha}(0)}$ and $\|\pmb{d}\|_{C^{l-2,\alpha}(0)}$.
\end{theorem}
~\\

\begin{theorem}[\textbf{$C^{k,\alpha}$ regularity}]\label{t-Cka-ell}
Let $0<\alpha<1$ and $u$ be a weak solution of \Cref{e.div-0}. Suppose that
\begin{equation*}
  \begin{aligned}
&\pmb{a}\in C^{k-1,\alpha}(0),~~\pmb{b}\in C^{k-2,\alpha}(0),
~~c\in C^{k-2, \alpha}(0),~~\pmb{d}\in C^{k-1,\alpha}(0),\\
&f \in C^{k-2,\alpha}(0),~~\pmb{f}\in C^{k-1,\alpha}(0).
  \end{aligned}
\end{equation*}

Then $u\in C^{k,\alpha}(0)$, i.e., there exists $P\in \mathcal{P}_k$ such that
\begin{equation*}\label{e.Cka.esti-ell}
  \|u-P\|^*_{L^2(B_r)}\leq C r^{k+\alpha}\left(\|u\|^*_{L^{2}(B_1)}+\|f\|_{C^{k-2,\alpha}(0)}
  +\|\pmb{f}\|_{C^{k-1,\alpha}(0)}\right), ~~\forall ~0<r<1,
\end{equation*}
\begin{equation*}
  \mathbf{\Pi}_{k-2} \left(P_{t}-(P^{a^{ij}}P_{i}+P^{d^j}P)_j+P^{b^i} P_{i}+P^cP-P^{f}+P^{f^i}_i\right)=0
\end{equation*}
and
\begin{equation*}
  |P(0)|+|DP(0)|+\cdots+|D^k P(0)|\leq C\left(\|u\|^*_{L^{2}(B_1)}+\|f\|_{C^{k-2,\alpha}(0)}
  +\|\pmb{f}\|_{C^{k-1,\alpha}(0)}\right),
\end{equation*}
where $C$ depends only on $n,\lambda, \Lambda, \alpha$, $\|\pmb{a}\|_{C^{k-1,\alpha}(0)}$, $\|\pmb{b}\|_{C^{k-2,\alpha}(0)}$, $\|c\|_{C^{k-2,\alpha}(0)}$ and $\|\pmb{d}\|_{C^{k-1,\alpha}(0)}$.
\end{theorem}

\begin{remark}\label{re6.1}
The proofs of above theorems are much simpler than that of parabolic equations. The energy inequality (Caccioppoli inequality) can provide the compactness directly by the Sobolev embedding and the  Fr\'{e}chet–Kolmogorov theorem (Kolmogorov–Riesz theorem) is not needed.

Certainly, in the proofs of above theorems, we will use $B_r$ instead of $Q_r$ and a polynomial $P\in\mathcal{P}_k$ means
\begin{equation*}
P(x)=\sum_{|\sigma|\leq k}\frac{a_{\sigma}}{\sigma!}x^{\sigma},~\sigma\in \mathbb{N}^n.
\end{equation*}
\end{remark}
~\\

Similar to the parabolic case, we have the following characterization of the nodal sets of solutions. Define
\begin{equation*}
  \mathcal{L}_k(u)=\left\{x\in \bar{B}_{1/2}:~~u(x)=\cdots=|D^{k-1}(x)|=0,~~D^ku(x)\neq 0 \right\}.
\end{equation*}
Then we have
\begin{theorem}\label{th1.4}
Let $l\geq 1$, $0<\alpha<1$ and $u$ be a weak solution of \Cref{e.div-0}. Suppose that
\begin{equation*}
  \begin{aligned}
    &\pmb{a}\in C^{l-1,\alpha}(\bar{B}_1),~~\pmb{b}\in C^{l-2,\alpha}(\bar{B}_1),~~c\in C^{l-3,\alpha}(\bar{B}_1),~~\pmb{d}\in C^{l-2,\alpha}(\bar{B}_1).
  \end{aligned}
\end{equation*}
Then
\begin{equation*}
  \begin{aligned}
\mathcal{L}_1(u)=\bigcup_{j=0}^{n-1} \mathcal{L}_1^j, \quad
\mathcal{L}_k(u)=\bigcup_{j=0}^{n-2} \mathcal{L}_k^j~~(k\geq 2),\\
  \end{aligned}
\end{equation*}
where each $\mathcal{L}_k^j$ is on a finite union of $j$-dimensional $C^{l,\alpha}$ manifolds.
\end{theorem}

\section{Appendix}\label{S.A}
In this appendix, we show the equivalence between the classical $C^{k,\alpha}(\bar{Q}_1)$ and the space defined via the pointwise $C^{k,\alpha}$ smoothness (see \Cref{d-f1}). Based on this equivalence, \Cref{co1.0} and \Cref{co1.1} follows easily.

First, we recall the classical $C^{k,\alpha}$ spaces in the parabolic distance (see \cite[P. 7-8]{MR0241822} and \cite[P. 46, Chapter IV.1]{MR1465184}).
\begin{definition}\label{de6.1}
Let $\Omega\subset \mathbb{R}^{n+1}$ be a bounded domain and $f:\Omega\rightarrow \mathbb{R}$. If $f\in C(\bar{\Omega})$, define
\begin{equation*}
\|f\|_{C(\bar{\Omega})}=\sup_{X\in \Omega}|f(X)|.
\end{equation*}
For $0<\alpha\leq 1$, if
\begin{equation*}
[f]_{C^{\alpha}(\bar{\Omega})}
:=\sup_{\mathop{X,Y\in \Omega,}\limits_{X\neq Y}}\frac{|f(X)-f(Y)|}{|X-Y|^{\alpha}}<\infty,
\end{equation*}
we say that $f\in C^{\alpha}(\bar{\Omega})$ ($f\in C^{0,1}(\bar{\Omega})$ for $\alpha=1$). We also introduce following semi-norm:
\begin{equation*}
[f]^{t}_{C^{\alpha}(\bar{\Omega})}:=\sup_{\mathop{(x,t),(x,s)\in \Omega,}\limits_{t\neq s}}
\frac{|f(x,t)-f(x,s)|}{|t-s|^{\alpha}}.
\end{equation*}

For $k\geq 1$, we say that $f\in C^k(\Omega)$ if $f$ has all derivatives of order $\leq k$ continuous in $\Omega$. Define
\begin{equation}\label{e6.1}
  [f]_{C^{k}(\bar{\Omega})}=\|D^kf\|_{C(\bar{\Omega})}+[D^{k-1}f]^{t}_{C^{1/2}(\bar{\Omega})},
\end{equation}
We say that $f\in C^{k}(\bar\Omega)$ if
\begin{equation*}
\|f\|_{C^{k}(\bar{\Omega})}:=\sum_{i=0}^{k-1}\|D^if\|_{C(\bar{\Omega})}+[f]_{C^k(\bar{\Omega})}<\infty.
\end{equation*}

Similarly, for $k\geq 1$ and $0<\alpha\leq 1$, define
\begin{equation}\label{e6.2}
  [f]_{C^{k,\alpha}(\bar{\Omega})}=
[D^kf]_{C^{\alpha}(\bar{\Omega})}+[D^{k-1}f]^{t}_{C^{(1+\alpha)/2}(\bar{\Omega})}.
\end{equation}
We say that $f\in C^{k,\alpha}(\bar\Omega)$ if
\begin{equation*}
\|f\|_{C^{k,\alpha}(\bar{\Omega})}:=\sum_{i=0}^{k}\|D^if\|_{C(\bar{\Omega})}
+[f]_{C^{k,\alpha}(\bar{\Omega})}<\infty.
\end{equation*}

If $f\in C^k(\Omega)$ and $X_0\in \Omega$. We define the parabolic $k$-th order Taylor polynomial of $f$ at $X_0$ as follows:
\begin{equation}\label{eA.10}
P^f_{X_0}(X)=\sum_{|\sigma|\leq k}\frac{1}{\sigma!}D^{\sigma}f(X_0)(X-X_0)^{\sigma},~\forall ~X\in \mathbb{R}^{n+1}.
\end{equation}
Note that
\begin{equation*}
D^{\sigma}P^f_{X_0}(X_0)=D^{\sigma}f(X_0), ~\forall ~|\sigma|\leq k.
\end{equation*}
\end{definition}

\begin{remark}\label{re.6.1}
Note that the definitions in \cref{e6.1} and \cref{e6.2} are different from the classical definitions in the usual Euclidean distance. If $k$ is even, the last terms in \cref{e6.1} and \cref{e6.2} can be removed.

Indeed, if $k$ is even, $[D^kf]_{C^{\alpha}}$ contains $k+\alpha$ order derivatives with respect to $x$ and $(k+\alpha)/2$ order derivatives with respect to $t$. However, if $k$ is odd, $[D^kf]_{C^{\alpha}}$ only contains $(k-1+\alpha)/2$ order derivatives with respect to $t$. Hence, we need an additional term $[D^{k-1}f]^{t}_{C^{(1+\alpha)/2}}$, which contains $(k-1)/2+(1+\alpha)/2=(k+\alpha)/2$ order derivatives with respect to $t$.
\end{remark}

With above notation, we have the following Taylor formula for smooth functions in parabolic distance.
\begin{lemma}\label{leA.3}
Suppose that $f\in C^k(Q_1)$. Then for any $(x_0,t_0)\in Q_1$,
\begin{equation*}
  \begin{aligned}
f&(x,t)-P^f_{x_0,t_0}(x,t)\\
=&\sum_{|\chi|=k}
\frac{1}{\chi!}D^{\chi}f(x_0+\theta_{\chi}(x-x_0),t)(x-x_0)^{\chi}\\
&+\sum_{\mathop{|\sigma|=k}\limits_{\sigma_{n+1}\geq 1}}
\frac{1}{\sigma!}D^{\sigma}f(x_0,t_0+\theta_{\sigma}(t-t_0))(x-x_0,t-t_0)^{\sigma}\\
 &+\sum_{|\sigma|=k-1}\frac{1}{\sigma!}
 \big(D^{\sigma}f(x_0,t_0+\theta_{\sigma}(t-t_0))-D^{\sigma}f(x_0,t_0)\big)(x-x_0,t-t_0)^{\sigma},
 ~\forall ~x\in Q_1,
  \end{aligned}
  \end{equation*}
where $P^f_{x_0,t_0}$ denotes the $(k-1)$-th order Taylor polynomial of $f$ at $(x_0, t_0)$ and $0\leq \theta_{\chi},\theta_{\sigma}\leq 1$.
\end{lemma}
\proof Without loss of generality, we assume $(x_0,t_0)=0$ and we need to prove
\begin{equation}\label{eA.9}
  \begin{aligned}
f(x,t)-P^f_0(x,t)=&\sum_{|\chi|=k}
\frac{1}{\chi!}D^{\chi}f(\theta_{\chi}x,t)x^{\chi}
+\sum_{\mathop{|\sigma|=k}\limits_{\sigma_{n+1}\geq 1}}
\frac{1}{\sigma!}D^{\sigma}f(0,\theta_{\sigma}t)(x,t)^{\sigma}\\
 &+\sum_{|\sigma|=k-1}
 \frac{1}{\sigma!}\left(D^{\sigma}f(0,\theta_{\sigma}t)-D^{\sigma}f(0)\right)(x,t)^{\sigma},~\forall ~(x,t)\in Q'_1.
  \end{aligned}
  \end{equation}
By the classical Taylor formula with respect to $x$, for any $(x,t)\in Q_1$,
\begin{equation*}
  \begin{aligned}
f(x,t)-P^f_{0,t}(x,t)=f(x,t)-\sum_{|\chi|\leq k-1} \frac{1}{\chi!}D^{\chi}f(0,t)x^{\chi}
=\sum_{|\chi|= k} \frac{1}{\chi!}D^{\chi}f(\theta_1x,t)x^{\chi},
  \end{aligned}
\end{equation*}
where $0\leq \theta_1\leq 1$. Similarly, by the Taylor formula for $D^{\chi}f(0,t)$ with respect to $t$, we have
\begin{equation*}
  \begin{aligned}
P^f_{0,t}&(x,t)-P_{0}(x,t)\\
= &\sum_{|\chi|\leq k-1} \frac{1}{\chi!}D^{\chi}f(0,t)x^{\chi}
-\sum_{|\chi|+2\tau\leq k-1}\frac{1}{\chi!\tau!}D^{\chi}D^{\tau}f(0)x^{\chi}t^{\tau}\\
=&\sum_{|\chi|\leq k-1} \frac{1}{\chi!}D^{\chi}f(0,t)x^{\chi}
-\sum_{|\chi|\leq k-1}\sum_{\tau\leq (k-1-|\chi|)/2} \frac{1}{\tau!}
 \left(\frac{1}{\chi!} D^{\tau}D^{\chi}f(0)x^{\chi}\right)t^{\tau}\\
=&\sum_{\mathop{|\sigma|=k}\limits_{\sigma_{n+1}\geq 1}}
\frac{1}{\sigma!}D^{\sigma}f(0,\theta_{\sigma}t)(x,t)^{\sigma}
+\sum_{|\sigma|=k-1}
 \frac{1}{\sigma!}\left(D^{\sigma}f(0,\theta_{\sigma}t)-D^{\sigma}f(0)\right)(x,t)^{\sigma},
  \end{aligned}
\end{equation*}
where $0\leq \theta_{\sigma}\leq 1$. By combining above estimates together, we obtain \cref{eA.9}.~\qed~\\

Next, we introduce following equivalence of norms for polynomials. Recall that $Q'(x_0,t_0,r):=B_r(x_0)\times (t_0-r^2,t_0+r^2)$.
\begin{lemma}\label{le1.1-2}
For $1\leq p,p'\leq +\infty$ and $k\geq 0$, we have:\\
(i) For any $X_0\in \mathbb{R}^{n+1}, r>0,~P\in \mathcal{P}_k$,
\begin{equation}\label{e.a.1}
\begin{aligned}
\|P\|^*_{L^{p,p'}(Q_r(X_0))}\simeq& \|P\|^*_{L^{p,p'}(Q'_r(X_0))}\simeq \sum_{i=0}^{k} r^i\|D^iP\|_{L^{\infty}(Q_r(X_0))}\simeq\sum_{i=0}^{k} r^i|D^iP(X_0)|\\
\end{aligned}
\end{equation}
(ii) For any $X_0\in Q_1, 0<r<1,~P\in \mathcal{P}_k,$
\begin{equation}\label{e.a.2}
\begin{aligned}
\|P\|^*_{L^{p,p'}(Q'_r(X_0))}\simeq& \|P\|^*_{L^{p,p'}(Q'_r(X_0)\cap Q_1)}
\simeq \sum_{i=0}^{k} r^i\|D^iP\|_{L^{\infty}(Q'_r(X_0)\cap Q_1)}.
\end{aligned}
\end{equation}
In above equations, $A\simeq B$ means
\begin{equation*}
C^{-1}A\leq B\leq CA,
\end{equation*}
where $C$ is a constant depending only on $n,k,p$ and $p'$.
\end{lemma}
\proof Clearly, $\mathcal{P}_{k}$ is a finite-dimensional vector space depending only on $n$ and $k$. Hence, any two norms in $\mathcal{P}_{k}$ are equivalent. In particular,
\begin{equation*}
\begin{aligned}
\|P\|^*_{L^{p,p'}(Q_1)}\simeq \|P\|^*_{L^{p,p'}(Q'_1)}\simeq \sum_{i=0}^{k} \|D^iP\|_{L^{\infty}(Q_1)}\simeq\sum_{i=0}^{k} |D^iP(0)|,
~\forall~P\in \mathcal{P}_k.
\end{aligned}
\end{equation*}
Note that the first norm and the second norm depend only on $n,p$ and $p'$; the third norm and the forth norm depend only on $n$ and $k$. Therefore, the equivalence constant $C$ depends only on $n,k,p$ and $p'$.

For a general $Q_r(X_0)$, by dilation and translation, we can transform $Q_r(X_0)$ to $Q_1$. Then \cref{e.a.1} holds with the same $C$ as in the case $Q_1$.

For \cref{e.a.2}, note that
\begin{equation*}
|Q'_r(X_0)\cap Q_1|\geq \frac{1}{2}|Q'_r(X_0)|,~\forall ~X_0\in Q_1,0<r<1.
\end{equation*}
Then \cref{e.a.2} can be proved similarly and we omit it.~\qed~\\

Now, we prove the equivalence in $Q_1$.
\begin{theorem}\label{a-l2.1}
Let $f\in L^{p,p'}(Q_1)$ ($1\leq p,p'\leq +\infty$). Then $f\in C^{k,\alpha}(\bar{Q}_1)(k\geq 0, 0<\alpha\leq 1)$ if and only if there exists a constant $C_0$ such that the following holds:

For any $X_0\in Q_{1}$, there exists $P_{X_0}\in \mathcal{P}_k$ such that
\begin{equation}\label{a-e3.1-1}
\|f-P_{X_0}\|^*_{L^{p,p'}(Q'(X_0,r)\cap Q_1)}\leq C_0r^{k+\alpha}, ~\forall ~0<r<1
\end{equation}
and
\begin{equation*}
  \begin{aligned}
K_0:=\sup_{X_0\in Q_1}  \sum_{i=0}^{k}|D^{i}P_{X_0}(X_0)|
+\sup_{\mathop{X_0\in Q_1,}\limits_{0<r<1}} r^{-(k+\alpha)}\|f-P_{X_0}\|^*_{L^{p,p'}(Q'(X_0,r)\cap Q_1)}<+\infty.
  \end{aligned}
\end{equation*}

Moreover, we have
\begin{equation}\label{a-e1.1-2}
C^{-1} \|f\|_{C^{k,\alpha}(\bar{Q}_1)}\leq K_0\leq C\|f\|_{C^{k,\alpha}(\bar{Q}_1)},
\end{equation}
where $C$ depends only on $n,k,p,p'$ and $\alpha$.
\end{theorem}
\proof We first prove the ``only if'' part. For any $X_0=(x_0,t_0)\in Q_1$, choose $P_{X_0}$ to be the $k$-th Taylor polynomial of $f$ at $X_0$ (cf. \cref{eA.10}). By applying the Taylor formula at $X_0$ (\Cref{leA.3}), we have
\begin{equation*}
  \begin{aligned}
|f(x&,t)-P_{X_0}(x,t)|\\
=&\big|\sum_{|\chi|= k} \frac{1}{\chi!}D^{\chi}f(\xi_{\chi})(x-x_0)^{\chi}-\sum_{|\chi|= k} \frac{1}{\chi!}D^{\chi}f(x_0,t_0)(x-x_0)^{\chi}\\
&+\sum_{\mathop{|\sigma|=k}\limits_{\sigma_{n+1}\geq 1}}
\frac{1}{\sigma!}D^{\sigma}f(\xi_{\sigma})(x-x_0,t-t_0)^{\sigma}
-\sum_{\mathop{|\sigma|=k}\limits_{\sigma_{n+1}\geq 1}}
\frac{1}{\sigma!}D^{\sigma}f(x_0,t_0)(x-x_0,t-t_0)^{\sigma}\\
 &+\sum_{|\sigma|=k-1}
 \frac{1}{\sigma!}\left(D^{\sigma}f(\xi_{\sigma})-D^{\sigma}f(x_0,t_0)\right)(x-x_0,t-t_0)^{\sigma}\big|\\
 \leq& C\|f\|_{C^{k,\alpha}(\bar{Q}_1)}\big(|(x,t)-(x_0,t_0)|^{\alpha}|x-x_0|^k
 +|(x,t)-(x_0,t_0)|^{k}|t-t_0|^{\alpha/2}\big)\\
 &+C\|f\|_{C^{k,\alpha}(\bar{Q}_1)}|(x,t)-(x_0,t_0)|^{k-1}|t-t_0|^{(1+\alpha)/2}\\
 \leq&C\|f\|_{C^{k,\alpha}(\bar{Q}_1)}|(x,t)-(x_0,t_0)|^{k+\alpha},~\forall ~(x,t)\in Q_1,
  \end{aligned}
\end{equation*}
where $C$ depends only on $n$ and $k$ and
\begin{equation*}
\xi_{\chi}=(x_0+\theta_{\chi}(x-x_0),t), \quad \xi_{\sigma}=(x_0,t_0+\theta_{\sigma}(t-t_0)), \quad
0\leq \theta_{\chi},\theta_{\sigma}\leq 1.
\end{equation*}

Note that
\begin{equation}\label{a-e2.8}
|Q'(X,r)\cap Q_1|\geq C |Q'(X,r)|,~\forall ~X\in Q_1,0<r<1,
\end{equation}
where $C$ depends only on $n$. Hence, for any $0<r<1$,
\begin{equation*}
\|f-P_{X_0}\|^*_{L^{p,p'}(Q'(X_0,r)\cap Q_1)}\leq
\|f-P_{X_0}\|_{L^{\infty}(Q'(X_0,r)\cap Q_1)}\leq C\|f\|_{C^{k,\alpha}(\bar{Q}_1)}r^{k+\alpha}.
\end{equation*}
That is, \cref{a-e3.1-1} holds and $K_0\leq C\|f\|_{C^{k,\alpha}(\bar{Q}_1)}$.

In the following, we prove the ``if'' part. By \cref{a-e3.1-1} and the Lebesgue differentiation theorem,
\begin{equation*}
f(X)=P_{X}(X), ~\forall~ a.e. ~X\in Q_1.
\end{equation*}
Redefine $f$ in a null set such that above equality holds everywhere.

Now, we prove that $f$ is differentiable at any $X_0=(x_0,t_0)\in Q_1$ ($t_0<0$). For any $X=(x,t)\in Q_1$ with $r:=|X-X_0|<1/4$, by \cref{a-e3.1-1} and \cref{a-e2.8},
\begin{equation*}
\begin{aligned}
\|P_{X}-P_{X_0}\|^*_{L^{p,p'}(Q'(X_0,r)\cap Q_1)}
\leq& \|f-P_{X}\|^*_{L^{p,p'}(Q'(X,2r)\cap Q_1)}
+\|f-P_{X_0}\|^*_{L^{p,p'}(Q'(X_0,r)\cap Q_1)}\\
\leq& CK_0r^{k+\alpha}.
\end{aligned}
\end{equation*}
By the equivalence of norms (see \Cref{le1.1-2}),
\begin{equation}\label{a-e2.3}
\sum_{i=0}^{k}r^i\|D^iP_{X}-D^iP_{X_0}\|_{L^{\infty}(Q'(X_0,r))}\leq CK_0r^{k+\alpha}.
\end{equation}
Thus,
\begin{equation}\label{a-e2.2}
\begin{aligned}
|f(X)-P_{X_0}(X)|=&|P_{X}(X)-P_{X_0}(X)|\leq CK_0r^{k+\alpha}=CK_0|X-X_0|^{k+\alpha}.
\end{aligned}
\end{equation}
Hence, by setting $t=t_0$ in above inequality, $f$ is differentiable with respect to $x$ and
\begin{equation}\label{a-e2.4}
Df(X_0)=DP_{X_0}(X_0).
\end{equation}

Next, we prove that $f$ is second order differentiable at $X_0$. That is, $f$ is differentiable with respect to $t$ and $Df$ is differentiable with respect to $x$. By setting $x=x_0$ in \cref{a-e2.2}, the first assertion holds and
\begin{equation*}
D_tf(X_0)=D_tP_{X_0}(X_0).
\end{equation*}
By \cref{a-e2.3} and \cref{a-e2.4},
\begin{equation*}
\begin{aligned}
|Df(X)-DP_{X_0}(X)|=&|DP_{X}(X)-DP_{X_0}(X)|\leq CK_0r^{k-1+\alpha}=CK_0|X-X_0|^{k-1+\alpha}.
\end{aligned}
\end{equation*}
Hence, by setting $t=t_0$ in above inequality, $Df$ is differentiable with respect to $x$ and
\begin{equation*}
D^2_xf(X_0)=D^2_xP_{X_0}(X_0).
\end{equation*}
Therefore, $f$ is second order differentiable at $X_0$ and
\begin{equation*}
D^2f(X_0)=D^2P_{X_0}(X_0).
\end{equation*}
By similar arguments, we can prove that $f$ is $k$-th order differentiable in $Q_1$ and
\begin{equation}\label{a-e1.2-2}
D^{\sigma}f(X_0)=D^{\sigma}P_{X_0}(X_0),~\forall ~|\sigma|\leq k, ~\forall ~X_0\in Q_1, t_0<0.
\end{equation}

Finally, for any $X_0,X\in Q_1$. If $r:=|X-X_0|< 1/4$, by \cref{a-e2.3},
\begin{equation}\label{a-e1.3-2}
\begin{aligned}
  |D^kf(X)-D^{k}f(X_0)|=&|D^kP_{X}(X)-D^kP_{X_0}(X_0)|\\
  =&|D^kP_{X}(X_0)-D^kP_{X_0}(X_0)|\\
  \leq& CK_0r^{\alpha}=CK_0|X-X_0|^{\alpha}.
\end{aligned}
\end{equation}
In addition, by setting $x=x_0$ and considering the $k-1$ order derivatives in above inequality, we also have
\begin{equation}\label{a-e1.3-4}
\begin{aligned}
  |D^{k-1}f(x_0,t)-D^{k-1}f(x_0,t_0)|  \leq CK_0r^{1+\alpha}=CK_0|t-t_0|^{(1+\alpha)/2}.
\end{aligned}
\end{equation}

If $|X-X_0|\geq 1/4$,
\begin{equation}\label{a-e1.3-3}
\begin{aligned}
  |D^{k}f(X)-D^{k}f(X_0)|\leq |D^{k}f(X)|+|D^{k}f(X_0)| \leq CK_0\leq CK_0|X-X_0|^{\alpha}.
\end{aligned}
\end{equation}
For $k-1$ order derivatives,
\begin{equation}\label{a-e1.3-5}
  |D^{k-1}f(x_0,t)-D^{k-1}f(x_0,t_0)|\leq CK_0\leq CK_0|t-t_0)|^{(1+\alpha)/2}.
\end{equation}
Note that \crefrange{a-e1.2-2}{a-e1.3-5} imply the left inequality in \cref{a-e1.1-2}. Hence, the theorem is proved.\qed~\\

Now, we prove that the pointwise estimates imply the classical $C^{k,\alpha}$ continuity.
\begin{theorem}\label{thA.3}
Let $k\geq 0, 0<\alpha\leq 1$ and $f\in L^{p,p'}(Q_1)$ ($1\leq p,p'\leq +\infty$). Suppose that there exists a constant $C_0$ such that the following holds:

For any $X_0\in Q_{1/2}$, there exists $P_{X_0}\in \mathcal{P}_k$ such that
\begin{equation}\label{eA.11}
\|f-P_{X_0}\|^*_{L^{p,p'}(Q(X_0,r))}\leq C_0r^{k+\alpha}, ~\forall ~0<r<1/2
\end{equation}
and
\begin{equation*}
  \begin{aligned}
K_0:=&\sup_{X_0\in Q_{1/2}}  \sum_{i=0}^{k}|D^{i}P_{X_0}(X_0)|
+\sup_{\mathop{X_0\in Q_{1/2},}\limits_{0<r<1/2}} r^{-(k+\alpha)}\|f-P_{X_0}\|^*_{L^{p,p'}(Q(X_0,r))}<+\infty.
  \end{aligned}
\end{equation*}

Then $f\in C^{k,\alpha}(\bar{Q}_{1/2})$ and
\begin{equation*}
\|f\|_{C^{k,\alpha}(\bar{Q}_{1/2})}\leq CK_0,
\end{equation*}
where $C$ depends only on $n,k,p,p'$ and $\alpha$.
\end{theorem}
\proof We only need to show that \cref{eA.11} implies \cref{a-e3.1-1} (with $Q_1$ replaced by $Q_{1/2}$ and $0<r<1$ replaced by $0<r<1/4$). Indeed, for any $X_0=(x_0,t_0)\in Q_{1/2}$ and $0<r<1/4$, there exists $t_1=t_0+\theta r^2$ for some $0\leq \theta \leq 1$ such that $X_1=(x_0,t_1)\in Q_{1/2}$ and $Q'(X_0,r)\cap Q_{1/2}\subset Q(X_1,2r)$. By \Cref{le1.1-2} and \cref{eA.11},
\begin{equation*}
  \begin{aligned}
\|P_{X_0}-P_{X_1}\|^*_{L^{p,p'}(Q'(X_0,r))}
\leq&C\|P_{X_0}-P_{X_1}\|^*_{L^{p,p'}(Q(X_0,r))}\\ \leq&C\|f-P_{X_0}\|^*_{L^{p,p'}(Q(X_0,r))}+C\|f-P_{X_1}\|^*_{L^{p,p'}(Q(x_0,t_1,2r))}\\
\leq&Cr^{k+\alpha}.
  \end{aligned}
\end{equation*}
Thus,
\begin{equation}\label{eA.14}
\begin{aligned}
\|f&-P_{X_0}\|^*_{L^{p,p'}(Q'(X_0,r)\cap Q_{1/2})}\\
\leq&\|f-P_{X_1}\|^*_{L^{p,p'}(Q'(X_0,r)\cap Q_{1/2})}
+\|P_{X_0}-P_{X_1}\|^*_{L^{p,p'}(Q'(X_0,r)\cap Q_{1/2})}\\
\leq&C\|f-P_{X_1}\|^*_{L^{p,p'}(Q(X_1,2r))}
+C\|P_{X_0}-P_{X_1}\|^*_{L^{p,p'}(Q'(X_0,r))}\\
\leq& Cr^{k+\alpha}.
\end{aligned}
\end{equation}
~\qed~\\

To characterize the nodal sets, we need the Whitney extension theorems. The first is the usual Whitney extension theorem in Euclidean space (see \cite{MR1501735} \cite[P. 177, Theorem 4]{MR0290095}, \cite[Theorem 2.3.6]{MR1065993} or \cite[Lemma 7.10]{MR2962060}.
\begin{theorem}\label{thA.1}
Let $k\geq 1$, $0<\alpha\leq 1$, $E\subset \mathbb{R}^n$ be a closed set and $f: E\to \mathbb{R}$. Suppose that for any $x\in E$, there exists $P_{x}\in \mathcal{P}_k$ such that
\begin{equation*}
  \begin{aligned}
&(i)~P_x(x)=f(x),~\forall ~x\in E;\\
&(ii)~|D^iP_x(x)-D^iP_y(x)|\leq K|x-y|^{k+\alpha-i},~\forall ~x,y\in E,~\forall ~ 0\leq i\leq k.
  \end{aligned}
\end{equation*}
Then $f$ extends to a $C^{k,\alpha}$ function $\tilde{f}$ in $\mathbb{R}^n$, i.e.,
\begin{equation*}
D^i\tilde{f}(x)=D^iP_x(x),~\forall ~x\in E,~\forall ~0\leq i\leq k.
\end{equation*}
\end{theorem}

The next is the Whitney extension in the parabolic distance. The $C^1$ extension has been proved in \cite[Proposition 6.2]{MR3385162}.
\begin{theorem}\label{thA.2}
Let $k\geq 1$, $0<\alpha\leq 1$, $E\subset \mathbb{R}^{n+1}$ be a closed set and $f: E\to \mathbb{R}$. Suppose that for any $X=(x,t)\in E$, there exists $P_{X}\in \mathcal{P}_k$ such that
\begin{equation}\label{eA.0}
  \begin{aligned}
&(i)~P_{X}(X)=f(X),~\forall ~X\in E;\\
&(ii)~|D^iP_{X}(X)-D^iP_{Y}(X)|\leq K|X-Y|^{k+\alpha-i},~\forall ~X,Y\in E,~\forall ~ 0\leq i\leq k.
  \end{aligned}
\end{equation}
Then $f$ extends to a $C^{k,\alpha}$ function $\tilde{f}$ in $\mathbb{R}^{n+1}$, i.e.,
\begin{equation*}
D^i\tilde{f}(X)=D^iP_{X}(X),~\forall ~X\in E,~\forall ~0\leq i\leq k.
\end{equation*}
\end{theorem}

Before proving the theorem, we recall the following two lemmas. The first is the Whitney decomposition (see \cite[P. 342]{MR3385162} and \cite[Chapter VI.1]{MR0290095}). In the following, we shall use cubes whose sides are parallel to the axes, i.e., a cube $Q$ (with center $(x,t)$ and side length $r$) can be written as
\begin{equation*}
Q=[x-r,x+r]^n\times [t-r^2,t+r^2]\quad\mbox{for some}~~(x,t)\in \mathbb{R}^{n+1},~r>0.
\end{equation*}
In addition, we use $Q^*$ to denote the cube with the same center as $Q$ and side length $(1+1/8)r$.
\begin{lemma}\label{leA.1}
Let $E\subset \mathbb{R}^{n+1}$ be a closed set. Then there exist sequences of $X_j\in E$ and cubes $Q_j$ whose interiors are mutually disjoint, such that
\begin{equation}\label{eA.1}
  \begin{aligned}
&(i)~E^c=\bigcup_{j\geq 1}Q_j=\bigcup_{j\geq 1}Q^*_j;\\
&(ii)~ C^{-1} \mathrm{diam}(Q_j)\leq d(Q_j,E)\leq C\mathrm{diam}(Q_j),~\forall ~j\geq 1;\\
&(iii)~d(Q_j, E)=d(Q_j,X_j),~\forall ~j\geq 1;\\
&(iv)~C^{-1} |X-X_j|\leq \mathrm{diam}(Q^*_j)\leq C|X-X_j|,~\forall ~j\geq 1,~\forall ~X\in Q^*_j;\\
&(v)~|X_0-X_j|\leq C|X_0-X|,~\forall ~j\geq 1,~\forall ~X_0\in E, ~\forall ~X\in Q^*_j;\\
&(vi)~X \mbox{ is contained in at most } N \mbox{ cubes } Q^*_j,~\forall ~X\in E^c,
  \end{aligned}
\end{equation}
where all constants $C$ and $N$ in above inequalities depend only on $n$.
\end{lemma}
~\\

The other is a partition of unity (see \cite[P. 342]{MR3385162} and \cite[Chapter VI.1]{MR0290095}).
\begin{lemma}\label{leA.2}
Let $E$ and $Q_j$ be as in \Cref{leA.1}. Then there exist a sequence of functions $\varphi_j$ such that
\begin{equation}\label{eA.6}
  \begin{aligned}
&(i)~ \varphi_j\in C^{\infty}_c(Q^*_j),~0\leq \varphi_j\leq 1,~\forall ~j\geq 1;\\
&(ii)~\sum_{j=1}^{\infty} \varphi_j\equiv 1\quad\mbox{ in}~~E^c;\\
&(iii)~\|D^i\varphi_j\|_{L^{\infty}(Q^*_j)}\leq C \mathrm{diam}(Q_j)^{-i},~\forall ~i\geq 1,j\geq 1,
  \end{aligned}
\end{equation}
where $C$ depends only on $n$ and $i$.
\end{lemma}
~\\

Now, we can give the~\\
\noindent\textbf{Proof \Cref{thA.2}.} We follow the idea in \cite[Chapter VI]{MR0290095}. In the following proof, $C$ always denotes a constant depending only on $n$ and $k$. Note first that
\begin{equation}\label{eA.3}
\begin{aligned}
|D^iP_{X}(Z)-D^iP_{Y}(Z)|\leq CK\sum_{l=0}^{k-i}|X-Y|^{k+\alpha-i-l}|X-Z|^l,\\
~\forall ~X,Y\in E,~Z\in \mathbb{R}^{n+1},~\forall ~ 0\leq i\leq k.
\end{aligned}
\end{equation}
Indeed, by rewriting the polynomials at $X$ (cf. \Cref{leA.3}), for any $0\leq i\leq k$,
\begin{equation}\label{eA.15}
  \begin{aligned}
D^iP_{X}(Z)-D^iP_{Y}(Z)
=&\sum_{|\sigma|\leq k-i}\frac{1}{\sigma!} D^{\sigma}D^i\left(P_{X}-P_{Y}\right)(X)(Z-X)^{\sigma}.
  \end{aligned}
\end{equation}
By the assumption (ii) in \cref{eA.0}, we obtain \cref{eA.3}.

Define the extension function as follows:
\begin{equation*}
\tilde{f}(X)= \left\{
  \begin{aligned}
    &f(X) &&\quad\mbox{if}~~X\in E;\\
    &\sum_{j=1}^{\infty} P_{X_j}(X)\varphi_j(X)&&\quad\mbox{if}~~X\in E^c,\\
  \end{aligned}
  \right.
\end{equation*}
where $X_j$ are as in \Cref{leA.1}. In the following, we show that $\tilde{f}$ is the desired extension. Obviously, $\tilde{f}\in C^{\infty}(E^c)$ by noting (vi) in \cref{eA.1}. We also use $P_{X}$ to denote the $k$-th order Taylor polynomial of $\tilde{f}$ at $X$ for any $X\in E^c$.

First, we prove that for any $X_0\in E$,
\begin{equation}\label{eA.2}
|\tilde{f}(X)-P_{X_0}(X)|\leq CK|X-X_0|^{k+\alpha},~\forall ~ X\in \mathbb{R}^n.
\end{equation}
Indeed, if $X\in E$, by the assumption (ii) in \cref{eA.0} (with $i=0$), the above holds clearly. If $X\in E^c$, there exist at most $N$ cubes containing $X$. For simplicity, we relabel them as $Q_1,...Q_N$ (similarly hereinafter). With the aid of (ii) in \cref{eA.6} and \cref{eA.3},
\begin{equation*}
  \begin{aligned}
|\tilde{f}(X)-P_{X_0}(X)|&=|\sum_{j=1}^{N} P_{X_j}(X)\varphi_j(X)-P_{X_0}(X)|\\
    &=|\sum_{j=1}^{N} \left(P_{X_0}(X)-P_{X_j}(X)\right)\varphi_j(X)|\\
    &\leq CK\sum_{j=1}^{N}|\varphi_j(X)|\sum_{l=0}^{k}|X_0-X_j|^{k+\alpha-l}|X_0-X|^l.
  \end{aligned}
\end{equation*}
From (v) in \cref{eA.1}, we arrive at \cref{eA.2}.

Next, we show that for any $X_0\in E$,
\begin{equation}\label{eA.5}
|D^i\tilde{f}(X)-D^iP_{X_0}(X)|\leq CK|X-X_0|^{k+\alpha-i},~\forall ~ X\in E^c,~\forall ~1\leq i\leq k.
\end{equation}
The proof is similar to that of \cref{eA.2}. For any $X\in E^c$ and $|\sigma|=i$, with the aid of (ii), (iii) in \cref{eA.6} and \cref{eA.3},
\begin{equation*}
  \begin{aligned}
|D^{\sigma}&\tilde{f}(X)-D^{\sigma}P_{X_0}(X)|\\
=&\left|\sum_{j=1}^{N}\sum_{\beta\leq \sigma}\binom{\sigma}{\beta}
D^{\beta}P_{X_j}D^{\sigma-\beta}\varphi_j-D^{\sigma}P_{X_0}\right|\\
=&\left|\sum_{j=1}^{N}
\left(D^{\sigma}P_{X_j}-D^{\sigma}P_{X_0}\right)\varphi_j\right|
+\left|\sum_{j=1}^{N}\sum_{\beta< \sigma}\binom{\sigma}{\beta} D^{\beta}P_{X_j}D^{\sigma-\beta}\varphi_j\right|\\
=&\left|\sum_{j=1}^{N}
\left(D^{\sigma}P_{X_j}-D^{\sigma}P_{X_0}\right)\varphi_j\right|
+\left|\sum_{j=1}^{N}\sum_{\beta< \sigma}\binom{\sigma}{\beta} \left(D^{\beta}P_{X_j}(X)-D^{\beta}P_{X_1}(X)\right)D^{\sigma-\beta}\varphi_j\right|\\
:=&J_1+J_2,
  \end{aligned}
\end{equation*}
where we have used
\begin{equation*}
\sum_{j=1}^{N}\sum_{\beta< \sigma}\binom{\sigma}{\beta} D^{\beta}P_{X_1}(X)D^{\sigma-\beta}\varphi_j=\sum_{\beta< \sigma}\binom{\sigma}{\beta} D^{\beta}P_{X_1}(X)D^{\sigma-\beta}\left(\sum_{j=1}^{N}\varphi_j\right)=0.
\end{equation*}

For $J_1$, by \cref{eA.3} (similar to the proof of \cref{eA.2}),
\begin{equation*}
J_1\leq CK|X-X_0|^{k+\alpha-i}.
\end{equation*}
For $J_2$, by \cref{eA.3} again and (iii) in \cref{eA.6},
\begin{equation*}
J_2\leq CK\sum_{j=1}^{N}\mathrm{diam}(Q_j)^{-i+|\beta|}
\sum_{l=0}^{k-|\beta|}|X_1-X_j|^{k+\alpha-|\beta|-l}|X_1-X|^l.
\end{equation*}
By (ii) in \cref{eA.1}, the diameters of these $N$ cubes are mutually comparable. Hence, with the aid of (iv) in \cref{eA.1},
\begin{equation*}
\mathrm{diam}(Q_j)\geq C|X-X_j|, \quad |X_1-X_j|\leq |X_1-X|+|X-X_j|\leq C|X-X_j|.
\end{equation*}
Thus,
\begin{equation*}
J_2\leq CK|X-X_j|^{k+\alpha-i}\leq CK|X-X_0|^{k+\alpha-i}.
\end{equation*}
By combining the estimates of $J_1,J_2$ together, we have \cref{eA.5}.

Since
\begin{equation*}
D^i\tilde{f}(X)=D^iP_{X}(X),~\forall ~X\in E^c,~\forall ~0\leq i\leq k,
\end{equation*}
by \cref{eA.5},
\begin{equation*}
|D^iP_{X}(X)-D^iP_{X_0}(X)|\leq CK|X-X_0|^{k+\alpha-i},~\forall ~X\in E^c,X_0\in E,~\forall ~1\leq i\leq k.
\end{equation*}
Then similar to the proof of \cref{eA.3}, we have
\begin{equation}\label{eA.7}
\begin{aligned}
|D^iP_{X}(Z)-D^iP_{Y}(Z)|\leq CK\sum_{l=0}^{k-i}|X-Y|^{k+\alpha-i-l}|X-Z|^l,\\
~\forall ~X\in E^c,Y\in E,~Z\in \mathbb{R}^{n+1},~\forall ~ 0\leq i\leq k.
\end{aligned}
\end{equation}

Next, we prove
\begin{equation}\label{eA.4}
  \left\{
  \begin{aligned}
    &|D^i\tilde{f}(X)|\leq \tilde C,~\forall ~0\leq i\leq k;\\
    &|D^{i}\tilde{f}(X)|\leq CKd(X,E)^{k+\alpha-i},~\forall ~i\geq k+1,
  \end{aligned}
 \quad~\forall ~X\in E^c, \right.
\end{equation}
where $\tilde{C}$ depends only on $n,k,K$, $|X-X_1|$ and $\|P_{X_1}\|$. The first inequality follows immediately from \cref{eA.5}. For the second inequality, since $P_{X_j}\in \mathcal{P}_k$, we have for any $|\sigma|=i\geq k+1$
\begin{equation*}
  \begin{aligned}
|D^{\sigma}\tilde{f}(X)|
=&\left|\sum_{j=1}^{N}\sum_{\beta\leq \sigma}\binom{\sigma}{\beta}
D^{\beta}P_{X_j}D^{\sigma-\beta}\varphi_j\right|\\
=&\left|\sum_{j=1}^{N}\sum_{|\sigma-\beta|\geq i-k}\binom{\sigma}{\beta} D^{\beta}P_{X_j}D^{\sigma-\beta}\varphi_j\right|\\
=&\left|\sum_{j=1}^{N}\sum_{|\sigma-\beta|\geq i-k}\binom{\sigma}{\beta} \left(D^{\beta}P_{X_j}(X)-D^{\beta}P_{X_1}(X)\right)D^{\sigma-\beta}\varphi_j\right|.
  \end{aligned}
\end{equation*}
Then similar to the estimate of $J_2$, we have
\begin{equation*}
|D^{\sigma}\tilde{f}(X)|\leq CK|X-X_1|^{k+\alpha-i}\leq CKd(X,E)^{k+\alpha-i}.
\end{equation*}

Next, we will prove that for any $\tilde X_0\in E^c$,
\begin{equation}\label{eA.8}
|\tilde{f}(X)-P_{\tilde X_0}(X)|\leq CK|X-\tilde X_0|^{k+\alpha},~\forall ~ X\in \mathbb{R}^n.
\end{equation}
Without loss of generality, we assume $\tilde X_0=0\in Q_1$ for some cube $Q_1$. For any $X\in \mathbb{R}^n$, if $X\in Q_1^*$, we can apply the Taylor formula at $\tilde X_0$ (\Cref{leA.3}). Then with the aid of \cref{eA.4},
\begin{equation*}
  \begin{aligned}
  |\tilde{f}(X)-P_{\tilde X_0}(X)|
=&\big|\sum_{|\chi|=k+1}
\frac{1}{\chi!}D^{\chi}u(\theta_{\chi}x,t)x^{\chi}
+\sum_{\mathop{k+1\leq|\sigma|\leq k+2}\limits_{\sigma_{n+1}\geq 1}}
\frac{1}{\sigma!}D^{\sigma}u(0,\theta_{\sigma}t)(x,t)^{\sigma}\big|\\
 \leq&CKd(X,E)^{\alpha-1}|X|^{k+1}+CKd(X,E)^{\alpha-2}|X|^{k+2}\\
 \leq&CK|X|^{k+\alpha},
  \end{aligned}
\end{equation*}
where we have used
\begin{equation*}
d(X,E)\geq C\mathrm{diam}(Q_1^*)\geq C|X|.
\end{equation*}

If $X\notin Q_1^*$,
\begin{equation*}
|X-X_1|\leq |X-\tilde X_0|+|X_1-\tilde X_0|\leq C|X-\tilde X_0|.
\end{equation*}
With the aid of \cref{eA.2} and \cref{eA.7},
\begin{equation*}
  \begin{aligned}
|\tilde{f}(X)-P_{\tilde X_0}(X)|\leq& |\tilde{f}(X)-P_{X_1}(X)|+|P_{X_1}(X)-P_{\tilde X_0}(X)|\\
 \leq&CK|X-X_1|^{k+\alpha}+CK\sum_{l=0}^{k}|X_1-\tilde{X}_0|^{k+\alpha-l}|X-\tilde X_0|^l\\
 \leq&CK|X-\tilde{X}_0|^{k+\alpha}.
  \end{aligned}
\end{equation*}
Therefore, \cref{eA.8} holds. By \Cref{a-l2.1}, $\tilde{f}\in C^{k,\alpha}(Q'_r)$ for any $r>0$.
~\qed~\\

As a corollary, we have
\begin{corollary}\label{coA.1}
Let $k\geq 0, 0<\alpha\leq 1$ and $f\in L^{p,p'}(Q_1)$ ($1\leq p,p'\leq +\infty$). Let $E\subset Q_{1/2}$ be a compact set. Suppose that there exists a constant $K$ such that the following holds:\\
For any $X_0\in E$, there exists $P_{X_0}\in \mathcal{P}_k$ such that
\begin{equation}\label{eA.13}
\|f-P_{X_0}\|^*_{L^{p,p'}(Q(X_0,r))}\leq Kr^{k+\alpha}, ~\forall ~0<r<1/2
\end{equation}
and
\begin{equation*}
  \begin{aligned}
K_0:=&\sup_{X_0\in E}  \sum_{i=0}^{k}|D^{i}P_{X_0}(X_0)|
+\inf\left\{K: \cref{eA.13}\mbox{ holds with }K \right\}<+\infty.
  \end{aligned}
\end{equation*}

Then $f$ extends to a $C^{k,\alpha}$ function $\tilde{f}$ in $Q'_1$ and
\begin{equation}\label{eA.12}
\|\tilde{f}\|_{C^{k,\alpha}(\bar{Q'}_1)}\leq CK_0,
\end{equation}
where $C$ depends only on $n,k$ and $\alpha$.
\end{corollary}
\proof We need to show (i) and (ii) in \Cref{thA.2}. Clearly, (i) holds and we only need to show (ii). Given
$X,Y\in E$, if $r:=|X-Y|\geq 1/8$, (ii) holds clearly. If $r< 1/8$, there exists $Z\in Q_{1/2}$ such that
\begin{equation*}
Q(Z,r)\subset Q(X,4r)\cap Q(Y,4r).
\end{equation*}
Then
\begin{equation*}
\begin{aligned}
\|P_X-P_{Y}\|^*_{L^{p,p'}(Q(Z,r))}
\leq&\|f-P_{X}\|^*_{L^{p,p'}(Q(Z,r))}
+\|f-P_{Y}\|^*_{L^{p,p'}(Q(Z,r))}\\
\leq&\|f-P_{X}\|^*_{L^{p,p'}(Q(X,4r))}
+\|f-P_{Y}\|^*_{L^{p,p'}(Q(Z,4r))}\\
\leq& Cr^{k+\alpha}.
\end{aligned}
\end{equation*}
By the equivalence of norms for polynomials (see \Cref{le1.1-2}),
\begin{equation*}
|D^iP_{X}(Z)-D^iP_{Y}(Z)|\leq Cr^{k+\alpha-i}=|X-Y|^{k+\alpha-i},~\forall ~ 0\leq i\leq k.
\end{equation*}
By a further rearrangement (see \cref{eA.15}),
\begin{equation*}
|D^iP_{X}(X)-D^iP_{Y}(X)|\leq C|X-Y|^{k+\alpha-i},~\forall ~ 0\leq i\leq k.
\end{equation*}
That is, the assumption (ii) in \Cref{thA.2} is satisfied. Hence, $f$ extends to a $C^{k,\alpha}$ function in $\mathbb{R}^{n+1}$. Finally, the estimate \cref{eA.12} is inferred from the proof of \Cref{thA.2} since all the constants in the proof depend only on $n,k,K_0$. ~\qed~\\

%\bibliographystyle{amsplain}
%\bibliographystyle{elsarticle-num}
%\bibliography{PDE}
%% \linenu
\printbibliography

\end{document}